\documentclass[12pt]{article}

\RequirePackage[OT1]{fontenc}
\usepackage{natbib}
\RequirePackage[colorlinks,citecolor=blue,urlcolor=blue]{hyperref}
\usepackage[english]{babel}
\usepackage{amsthm}
\usepackage{amsmath}
\usepackage{amsfonts}
\usepackage{amssymb}
\usepackage{graphicx}
\usepackage{enumerate}
\usepackage{color}
\usepackage{array}
\usepackage{bbm}
\usepackage{multirow}
\usepackage{rotating}
\usepackage{xr}
\usepackage{setspace}

\theoremstyle{plain}
\newtheorem{theorem}{Theorem}

\theoremstyle{definition}

\begin{document}
\thispagestyle{empty}
\baselineskip=28pt
\vskip 5mm
\begin{center} {\Large{\bf On the Tail Behaviour of Aggregated Random Variables}}
\end{center}

\baselineskip=12pt
\vskip 5mm

\begin{center}
\large
Jordan Richards$^{12*}$ and Jonathan A. Tawn$^1$
\end{center}

\footnotetext[1]{
\baselineskip=10pt STOR-i Centre for Doctoral Training, Department of Mathematics and Statistics, Lancaster University, LA1 4YR, UK.}
\footnotetext[2]{
\baselineskip=10pt Computer, Electrical and Mathematical Sciences and Engineering Division,
King Abdullah University of Science and Technology (KAUST), Thuwal, 23955-6900, Saudi Arabia. \\$^*$E-mail: jordan.richards@kaust.edu.sa}

\baselineskip=17pt
\vskip 4mm
\centerline{\today}
\vskip 6mm

\begin{center}
{\large{\bf Abstract}}
\end{center}
In many areas of interest, modern risk assessment requires estimation of the extremal behaviour of sums of random variables. We derive the first order upper-tail behaviour of the weighted sum of bivariate random variables under weak assumptions on their marginal distributions and their copula. The extremal behaviour of the marginal variables is characterised by the generalised Pareto distribution and their extremal dependence through subclasses of the limiting representations of Ledford and Tawn \cite{Ledford1997} and Heffernan and Tawn \cite{Heffernan2004}. We find that the upper-tail behaviour of the aggregate is driven by different factors dependent on the signs of the marginal shape parameters; if they are both negative, the extremal behaviour of the aggregate is determined by both marginal shape parameters and the coefficient of asymptotic independence \citep{Ledford1996}; if they are both positive or have different signs, the upper-tail behaviour of the aggregate is given solely by the largest marginal shape. We also derive the aggregate upper-tail behaviour for some well known copulae which reveals further insight into the tail structure when the copula falls outside the conditions for the subclasses of the limiting dependence representations.
\baselineskip=16pt

\par\vfill\noindent
{\bf Keywords:} aggregation; coefficient of asymptotic independence; conditional extreme;   copula models; extreme value theory; generalised Pareto distribution.\\

\pagenumbering{arabic}
\baselineskip=24pt

\newpage


 \section{Introduction}
  \label{intro}
  The extremal behaviour of aggregated data is of importance in two key areas of risk management; financial portfolio optimisation and fluvial flooding. In financial risk management, it is standard practice to aggregate over returns from several assets in a portfolio in an attempt to mitigate investment risk. It is important that the uncertainty surrounding the tail behaviour of the aggregate is assessed so that the risk of large negative cumulative returns can be quantified \citep{Hauksson2001,Chen2019,Embrechts2015}. For flood risk management, consider that fluvial floods are typically caused by prolonged extreme precipitation over a catchment area; more succinctly, precipitation aggregated both spatially and temporally \citep{Coles1996,Eggert2015}. In both cases, the assumption of independence within the multivariate variable of interest is unlikely to hold. We derive the first order behaviour of the upper-tail of a weighted sum of a bivariate random vector with different marginal tail behaviours and extremal dependence structures and demonstrate that both factors have a significant effect on the extremal behaviour of the aggregate variable.  \par
We define the aggregate $R$ as a weighted sum of the components of a random vector ${\mathbf{X}=(X_1,\dots,X_d)}$, with marginal distribution functions $F_i$ for $i \in\{1,\dots,d\}$, as
\begin{equation}
\label{sum}
R=\sum^{d}_{i=1}\omega_i X_i,
\end{equation}
with weights $\boldsymbol{\omega}=\{\omega_i; 0 \leq \omega_i \leq 1, \sum^d_{i=1} \omega_i=1\},$ and where components of $\mathbf{X}$ are all positive and not necessarily independent and identically distributed and $\mathbf{X}$ has a joint density. Dependence between components can be described using copulae, Sklar's theorem \citep{Nelsen2006}. The joint distribution function of $\mathbf{X}$ can be uniquely written as $F(\mathbf{x})=C\{F_1(x_1),\dots,F_d(x_d)\}$, for $\mathbf{x} \in \mathbb{R}^d$, where $C$ is the copula, i.e., some multivariate distribution function $C :[0,1]^d \rightarrow [0,1]$ on uniform margins. Our interest lies in the tail behaviour of $R$, which we quantify by considering $\Pr\{R \geq r\}$ as $r \rightarrow r^F$, where $r^F \leq \infty$ is the upper-endpoint of $R$, and how this behaviour is driven by the marginal tails and dependence structure of $\mathbf{X}$.
Modelling the marginal tails of a random vector $\mathbf{X}$ has been widely studied, see \cite{Pickands1975,Davison1990} and \citet{coles2001}. The typical approach is to assume that there exists a threshold $u_i$ for each $X_i$, such that the distribution of $(X_i-u_i)| ( X_i > u_i)$ is characterised by a generalised Pareto distribution, denoted GPD$(\sigma_i,\xi_i),$ with distribution function
\begin{equation}
\label{GPDCDF}
H_i(x)=\begin{cases}
1-\left(1+\xi_i x/\sigma_i\right)_+^{-1/\xi_i},\;\;\;&\xi_i \neq 0,\\
1-\exp\left(-x/\sigma_i\right), &\xi_i = 0,
\end{cases}
\end{equation}
for $x>0$, scale parameter $\sigma_i > 0$, shape parameter $\xi_i \in \mathbb{R}$ and where ${z}_+=\max\{0,z\}$. The operator $z_+$ forces $X_i$ to have upper-endpoint $x_i^F =  u_i-\sigma_i/\xi_i$ if and only if $\xi_i \leq 0$ and the shape parameter $\xi_i$ controls the heaviness of the upper-tails of $X_i$: for $\xi_i >0,\; \xi_i = 0$ and $\xi_i < 0$, we have that $X_i$ has heavy, exponential and bounded, upper-tails, respectively. It is important to make the distinction between these three cases as we show that the sign of the marginal shape parameters, $\xi_i$, has a large effect on the tail behaviour of $R$. We focus on the bivariate sum $R=X_1+X_2$, where $X_i\sim \mbox{GPD}(\sigma_i,\xi_i)$ and $X_i>0$ for $i\in\{1,2\}$, i.e., setting $F_i = H_i$, and with some specified joint distribution on $(X_1,X_2)$; the choice of $u_i = 0$ for $i\in\{1,2\}$ is discussed in Section~\ref{motivation-sec}. \par
It remains to specify the dependence structure between $X_1$ and $X_2$ which leads to large $R$. The dependence between extreme values of variables is classified as either asymptotic dependence or asymptotic independence with respective measures of dependence: $\chi$ the coefficient of asymptotic dependence and $\bar{\chi}$ the coefficient of asymptotic independence \citep{coles1999}. The former is defined as
\begin{equation}
\label{chiEq1}
\chi=\lim_{q \uparrow 1} \Pr\{F_1(X_1)>q|F_2(X_2)>q\}.
\end{equation}
For $\chi = 0$ and $\chi > 0$, we have asymptotic independence and asymptotic dependence, respectively, with $\chi$ increasing with strength of extremal dependence. Conversely, \cite{Ledford1996} characterise asymptotic independence between $X_1$ and $X_2$ through the assumption that 
\begin{equation}
\label{etaEq1}
\Pr\left\{F_1(X_1) > 1-1/u, F_2(X_2) > 1-1/u\right\}= \mathcal{L}(u)u^{-1/\eta},
\end{equation}
where $0 < \eta \leq 1,\;\mathcal{L}(u)$ is a positive slowly varying function as $u\rightarrow \infty$ and $\bar{\chi}=2\eta - 1$, so $-1 < \bar{\chi} \leq 1$. If $\bar{\chi} = 1$ and $\mathcal{L}(u)$ tends to a positive constant as $u \rightarrow \infty$, we have asymptotic dependence, and for $\bar{\chi} \in [0,1)$ we have asymptotic independence with weakening strength of dependence as $\bar{\chi}$ decreases. We consider two special cases of these extremal dependence classes, namely perfect positive dependence with $\chi =1$ in \eqref{chiEq1} and $\eta=1$, and independence with $\chi =0$ and $\eta=1/2$. In both cases, $\mathcal{L}(u)=1$ for $u> 1$. \par
Previous studies on the tail behaviour of aggregated random variables focus on the effects of the marginal distributions, with limited cases of the dependence structure being considered. Numerous studies on the sum of independent $(\chi = 0,\;\bar{\chi} = 0)$ Pareto random variables, corresponding to GPD random variables with $\xi = 1$, have been conducted, see \cite{Ramsay2008,Nguyen2015}. The tail behaviour of weighted sums of Pareto random variables, where the weights are random and exhibit dependence is modelled using elliptical distributions, is studied by \cite{goovaerts2005tail} and  \cite{opitz2016modeling} describes the relationship between marginal exceedance probabilities for both an exponential-tailed Laplace random vector and its sum. The exact distribution for sums of independent exponential random variables with nonhomogenous, i.e., different, marginal scale parameters, is studied by \citet{Nadarajah2008} and \cite{NadarajahKotz2008,Nadarajah2018} extend this framework to independent GPD margins.  A further derviation of the distribution of $R$ with GPD margins and a Clayton copula $(\chi > 0,\;\bar{\chi} = 1)$, see \cite{Ghosh2020}, is provided by \cite{Nadarajah2006}. For asymptotically independent variables, \cite{mitra2009aggregation} study the behaviour of bivariate aggregates with exponential upper-tails, i.e., $\xi=0$. \par Under a general assumption that $\chi > 0$ and that the shape parameters are equal, studies that focus on the extremal behaviour of $R$ include \cite{Coles1994} and \cite{kluppelberg2008} and where $R$ is an integral of a stochastic process by \cite{Coles1996, ferreira2012exceedance} and \cite{Engelke2019}, with the latter studied numerically by \cite{richards2021modelling} for $\chi \geq 0$. Further extensions to asymptotically independent structures has been made by \cite{engelke2019extremal}, who study the relationship between the relative tail decay rates of the bivariate sum $R$ and random vector $(X_1/R,X_2/R)$, and the corresponding values of $\chi$ and $\eta$ for $(X_1,X_2)$; however, these are general results and do not link the marginal shapes to the tail decay rate of $R$. Other general results for the tail behaviour of sums include extensions of Breiman's lemma \citep{breiman1965some}, which link the decay rate of a multivariate regularly varying random vector to the decay rate of the sum of its components, see \cite{fougeres2012risk,li2018joint}. \par
There are important gaps in the literature for the tail behaviour of $R$ relating to unequal marginal shape parameters and copulae with $\chi = 0$ and $\bar{\chi} < 1$. The case where $\bar{\chi} < 0$ implies negative dependence between $X_1$ and $X_2$; this case is also absent from the literature, but we primarily constrain our focus to $\bar{\chi} \geq 0$.\par 
The paper is structured as follows. In Section~\ref{motivation-sec}, we conduct a preliminary investigation into the upper-tail behaviour of $R$ with results that motivate our modelling choices for $(X_1,X_2)$. Section~\ref{LimRes-sec} introduces our preliminary model set-up and the results that follow by modelling dependence in $(X_1,X_2)$ using limit models that are specified there which cover both asymptotic dependence and asymptotic indepedence cases; these results are easily interpretable and give a strong insight into the tail behaviour of the aggregate. In Section \ref{cop-sec}, we provide examples of our results for widely used copulae and give further insight into the tail behaviour of $R$ when the dependence in $(X_1,X_2)$ does not satisfy the conditions detailed in Section~\ref{Results-sec}. We apply our results to climate data in Section~\ref{Sec-Appl} and provide further discussion in Section~\ref{discuss-sec}. \ref{proof_sec} provides the proofs of the results in Section~\ref{Results-sec}; for full details see \cite{richardsthesis}.
\section{Motivation}
\label{motivation-sec}
We explore the upper-tail of $R$ numerically using Monte-Carlo methods for copulae with a range of $\chi$ and $\bar{\chi}$ values; this is to motivate the form in which we present the results in Section~\ref{Results-sec} and our choice of dependence model for $(X_1,X_2)$. We consider two copulae based on the bivariate extreme value copula, see \cite{Tawn1988} and \cite{GudendorfBook}. An example of a bivariate extreme value copula is the logistic model,
 \begin{equation}
 \label{logistic}
 C_L(u,v)=\exp \left\{-\left[(-\log u)^{1/\gamma}+(-\log v)^{1/\gamma}\right]^\gamma\right\},\;\;\;\;u,v \in [0,1],
 \end{equation}
 where $\gamma \in [0,1)$; where here we avoid the case $\gamma = 1$ which is the independence copula, but allow $\gamma = 0$, taken as the limit in \eqref{logistic} as $\gamma \rightarrow 0$. From \eqref{chiEq1} and \eqref{etaEq1}, this copula gives values $\chi = 2-2^\gamma>0$ and $\bar{\chi} = 1$, and the variables are asymptotically dependent with the strength of asymptotic dependence decreasing with $\gamma$ increasing. Inverting this copula gives the inverted-logistic copula which is asymptotically independent, see \cite{wadsworth2012}. This is defined through its survival copula, 
  \begin{equation}
 \label{invlogistic}
 \bar{C}_{IL}(u,v)=\exp \left\{-\left[\left(-\log(1-u)\right)^{1/\gamma}+\left(-\log(1-v)\right)^{1/\gamma}\right]^\gamma\right\},\;\;\;\;u,v \in [0,1],
 \end{equation}
 where $\gamma \in (0,1]$. In contrast to the logistic copula, we have $\chi = 0$ and $\bar{\chi} = 2^{1-\gamma}-1$, with strength of asymptotic independence increasing as $\gamma$ decreases. \par
  \begin{figure}[h!]
\centering
\begin{minipage}{0.44\linewidth} 
\centering
\includegraphics[width=\linewidth,angle=-90]{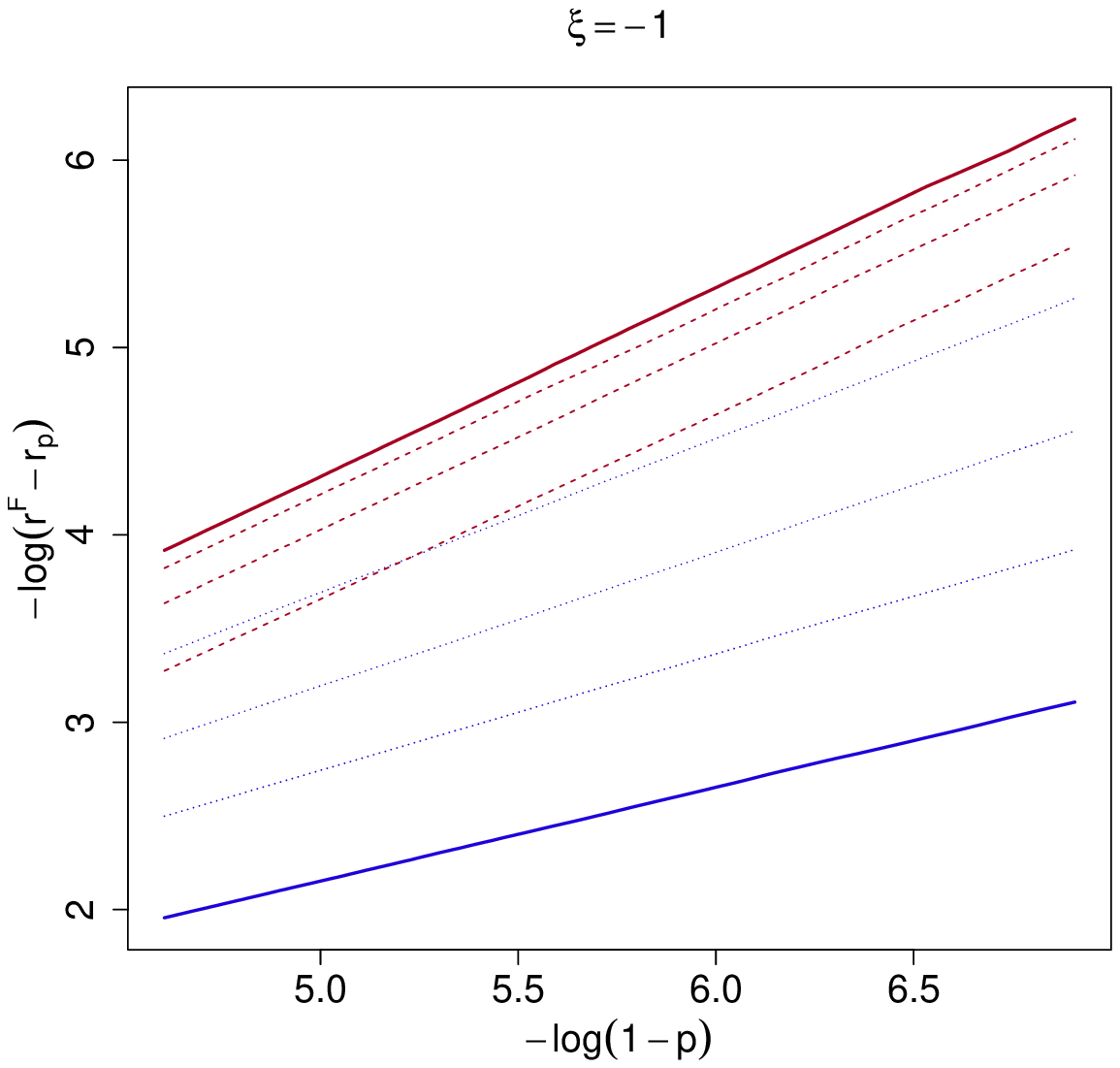} 
\end{minipage}
\begin{minipage}{0.44\linewidth} 
\centering
\includegraphics[width=\linewidth,angle=-90]{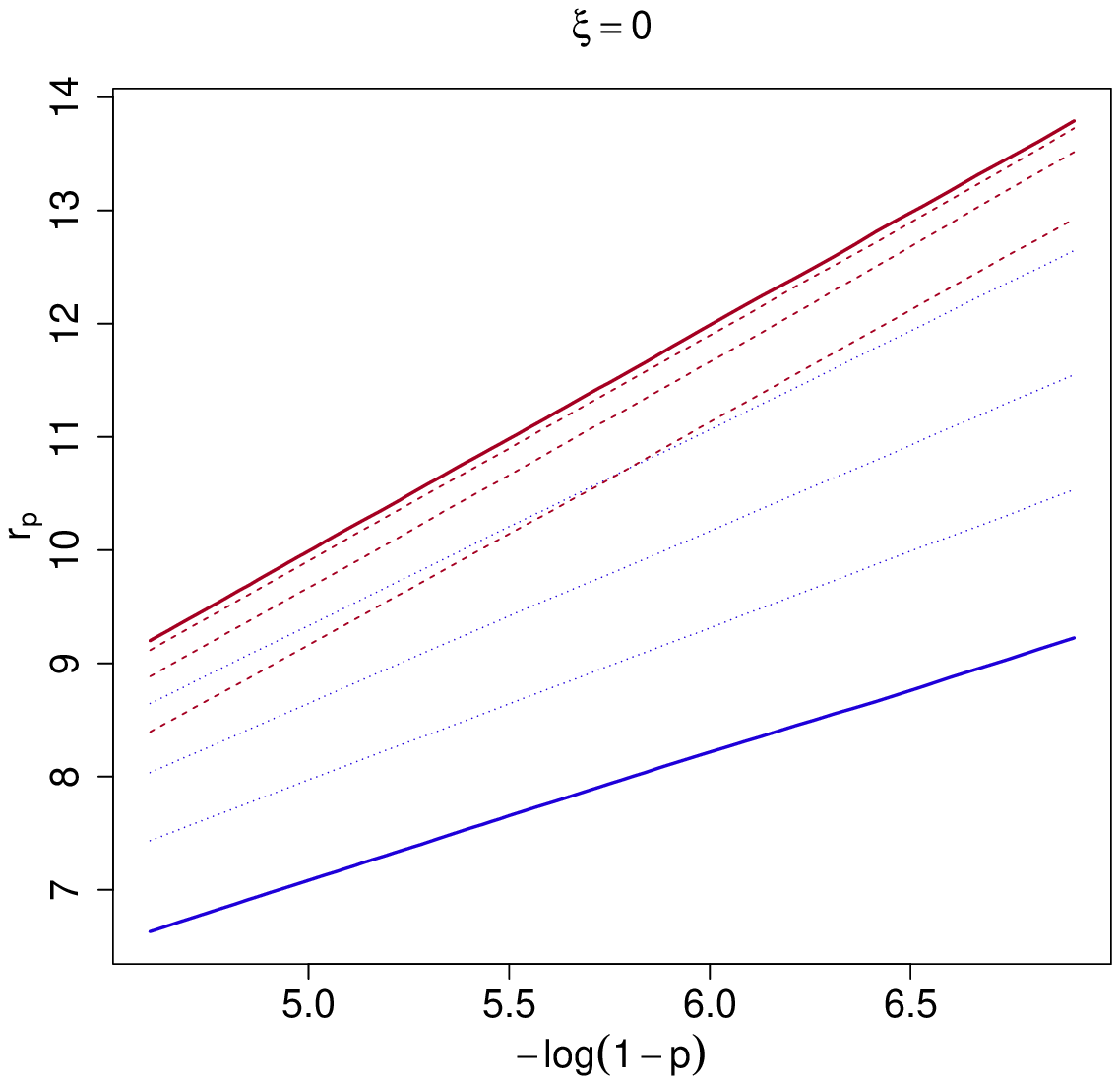}  
\end{minipage}
\begin{minipage}{0.44\linewidth} 
\centering
\includegraphics[width=\linewidth,angle=-90]{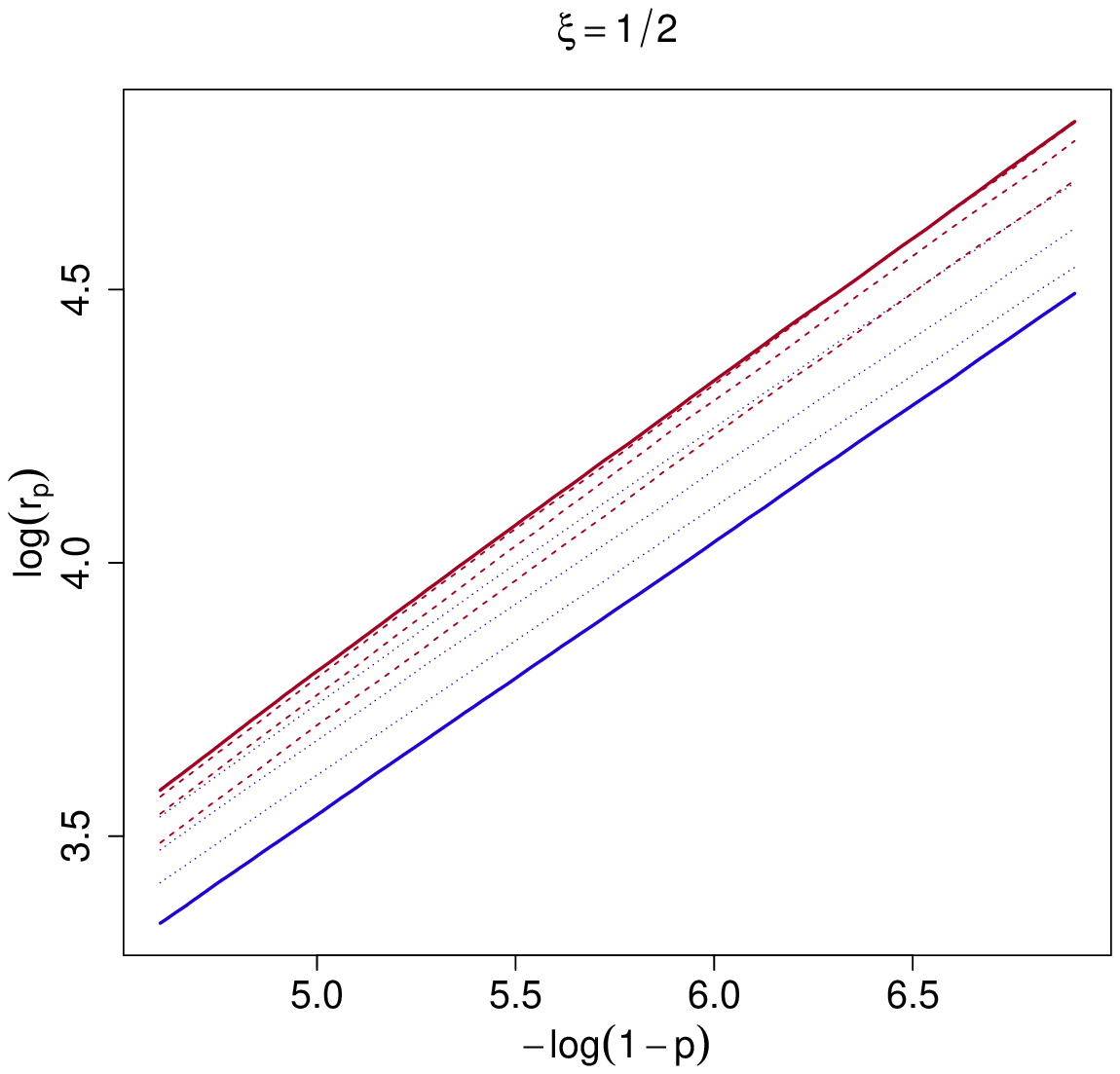}  
\end{minipage}
\begin{minipage}{0.44\linewidth} 
\centering
\includegraphics[width=\linewidth,angle=-90]{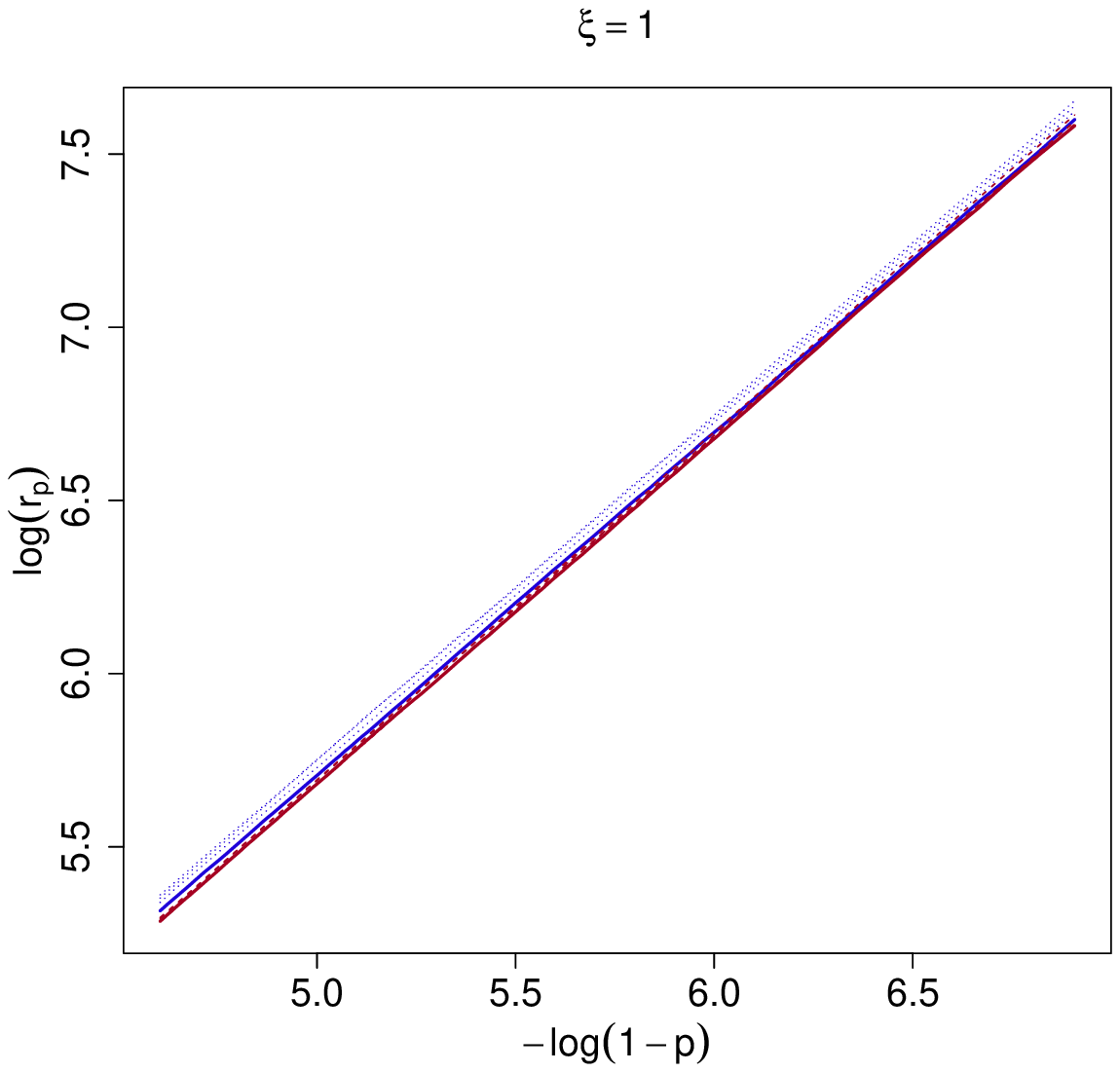} 
\end{minipage}
\caption{Quantiles $r_p$ of $R$; the sum of two GPD$(1,\xi)$ random variables, with copula \eqref{logistic} in red and \eqref{invlogistic} in blue and for $\xi = -1,0,1/2,1$ and $\gamma = 0.3,0.5,0.9$ and $p\in[0.99,0.999]$. To emphasise their similarities, these are displayed on the scales $-\log(r^F-r_p),\; r_p$ and $\log(r_p)$ for $\xi < 0,\; \xi = 0$ and $\xi > 0$ respectively, where $r^F$ is the upper-endpoint of $R$. Solid lines correspond to perfect dependence and independence, and the values on the $y-$axis decrease in each plot with increasing $\gamma$. Curves are estimated using Monte Carlo methods, with samples taken to be sufficiently large that any observed differences in the plot are statistically significant.}
\label{SimQuans}
\end{figure}
Fig.~\ref{SimQuans} provides simulated quantiles of samples of size $1\times 10^7$ for $R=X_1+X_2$, where $X_1,X_2 \sim \mbox{GPD}(1,\xi)$ with the bivariate extreme value logistic, and inverted logistic, copulae, defined in \eqref{logistic} and  \eqref{invlogistic}, respectively, for selected values of $\xi$ and the copulae parameter $\gamma$. Quantiles $r_p$, where $\tilde{F}_R(r_p) = p$ for $\tilde{F}_R$ the empirical distribution of $R$, are given for $p$ close to $1$. We observe that growth of the quantiles of $R$ is affected by both the underlying dependence in $(X_1,X_2)$ and the marginal shape parameters. The scales of the axes in Fig.~\ref{SimQuans} are chosen so that the approximate slope of the lines reveal the shape parameter of $R$. To illustrate this, we first define asymptotic notation ($\sim$); for any functions $f:\mathbb{R}\rightarrow \mathbb{R}$, $g:\mathbb{R}\rightarrow \mathbb{R}$ and fixed $a\in \mathbb{R}$, the relation $f(x)\sim g(x)$ as $x\rightarrow a$ holds if and only if $\lim_{x\rightarrow a}f(x)/g(x)=1$. Then when $\Pr\{R \geq r_p\} = 1-p$, we have
\begin{equation}
\label{eqRformquants}
-\log(1-p) \sim \begin{cases}
\xi_R^{-1}\log(r_p)-\log(K_1),\;\;&\text{if}\;\; \xi_R > 0,\\
\sigma_R^{-1}r_p-\log(K_2),\;\;&\text{if}\;\; \xi_R = 0,\\
-\xi_R^{-1}\log(r^F-r_p)-\log(K_3)-\xi_R^{-1}\log(r^F),\;\;&\text{if}\;\; \xi_R < 0,
\end{cases}
\end{equation}
as $p \rightarrow 1$ for constants $K_1,K_2,K_3 >0$. Thus, with the axes scaling used in Fig. \ref{SimQuans}, we expect the slope of each quantile curve to be approximately $1/\xi_R, \;1/\sigma_R$ and $-1/\xi_R$ if $\xi_R >0, \; \xi_R = 0$ and $\xi_R < 0$, respectively, for sufficiently large $R$. \par
Relationship \eqref{eqRformquants} and Fig.~\ref{SimQuans} reveal interesting preliminary insights into the upper-tail behaviour of $R$. For $\xi >0$, we find that the slopes in Fig.~\ref{SimQuans} are approximately equal; implying that the dependence structure has no significant effect on $\xi_R$. For $\xi \leq 0$, the reverse is true; for $\xi=0$, we observe that for the asymptotically independent copula \eqref{invlogistic}, then $\sigma_R$ changes with the strength of dependence; a similar property can be observed for $\xi < 0$, albeit given a change in $\xi_R$. In both cases, the slopes remain approximately equal for the quantiles derived using the asymptotically dependent copula \eqref{logistic}, which implies that some of the structure in $\xi_R$ is driven by the strength of asymptotic independence, rather than the degree of asymptotic dependence. We note that our empirical findings are in full agreement with the results for the theoretical upper-tail behaviour of $R$ that we detail in Section~\ref{Results-sec}. \par
\begin{figure}[h!]
\centering
\begin{minipage}{0.32\linewidth} 
\centering
\includegraphics[width=\linewidth,angle=-90]{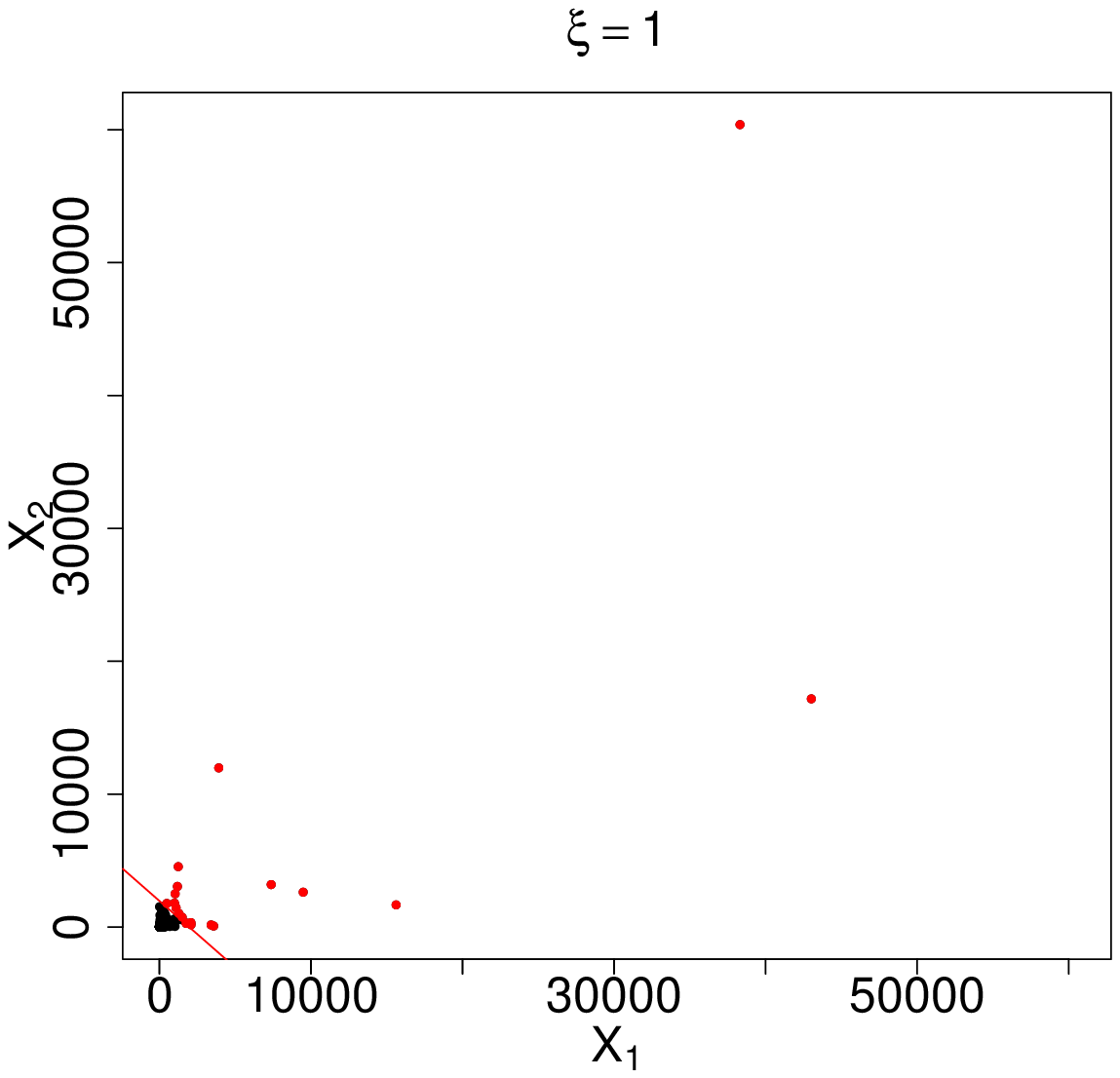} 
\end{minipage}
\begin{minipage}{0.32\linewidth} 
\centering
\includegraphics[width=\linewidth,angle=-90]{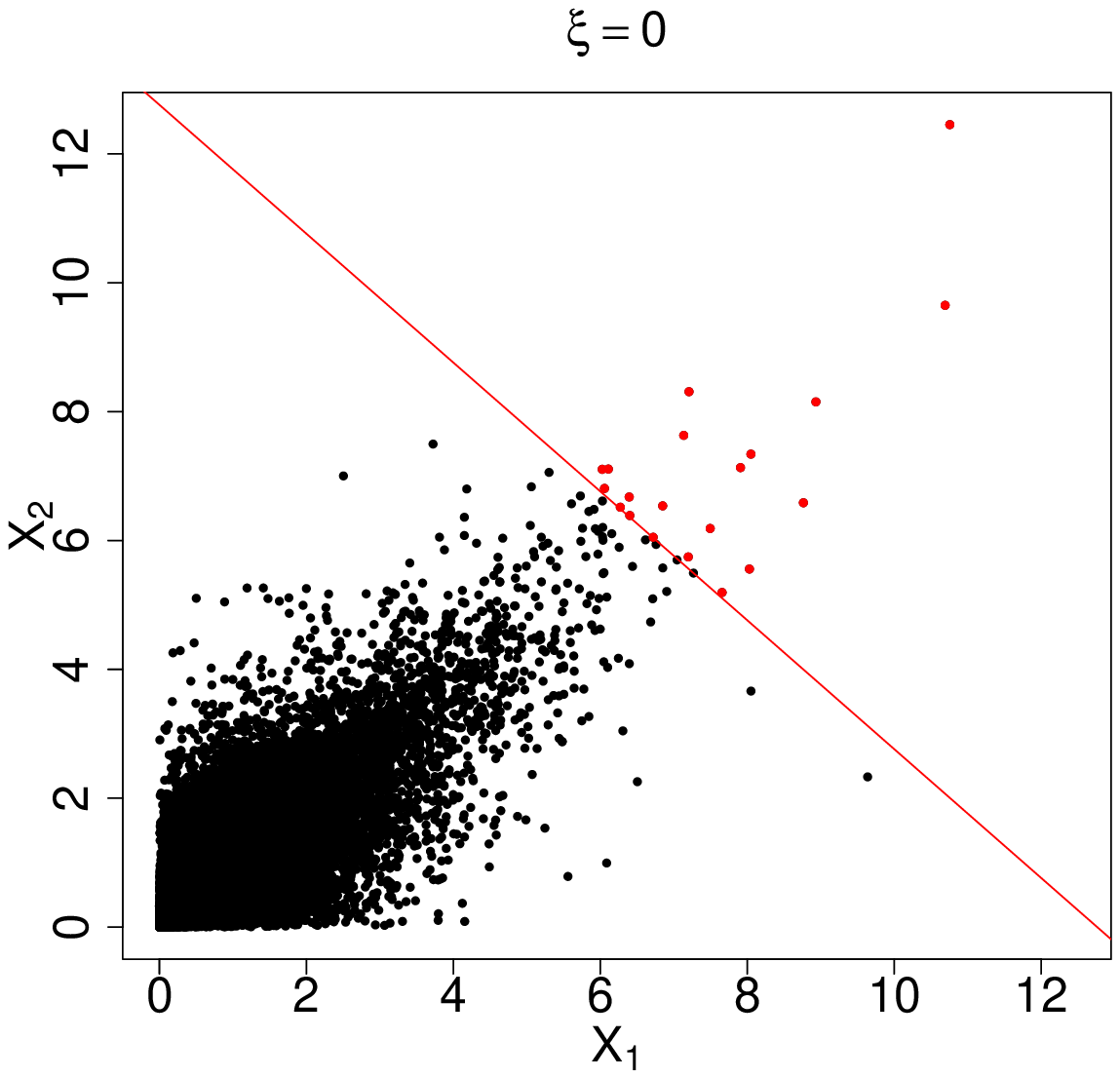}  
\end{minipage}
\begin{minipage}{0.32\linewidth} 
\centering
\includegraphics[width=\linewidth,angle=-90]{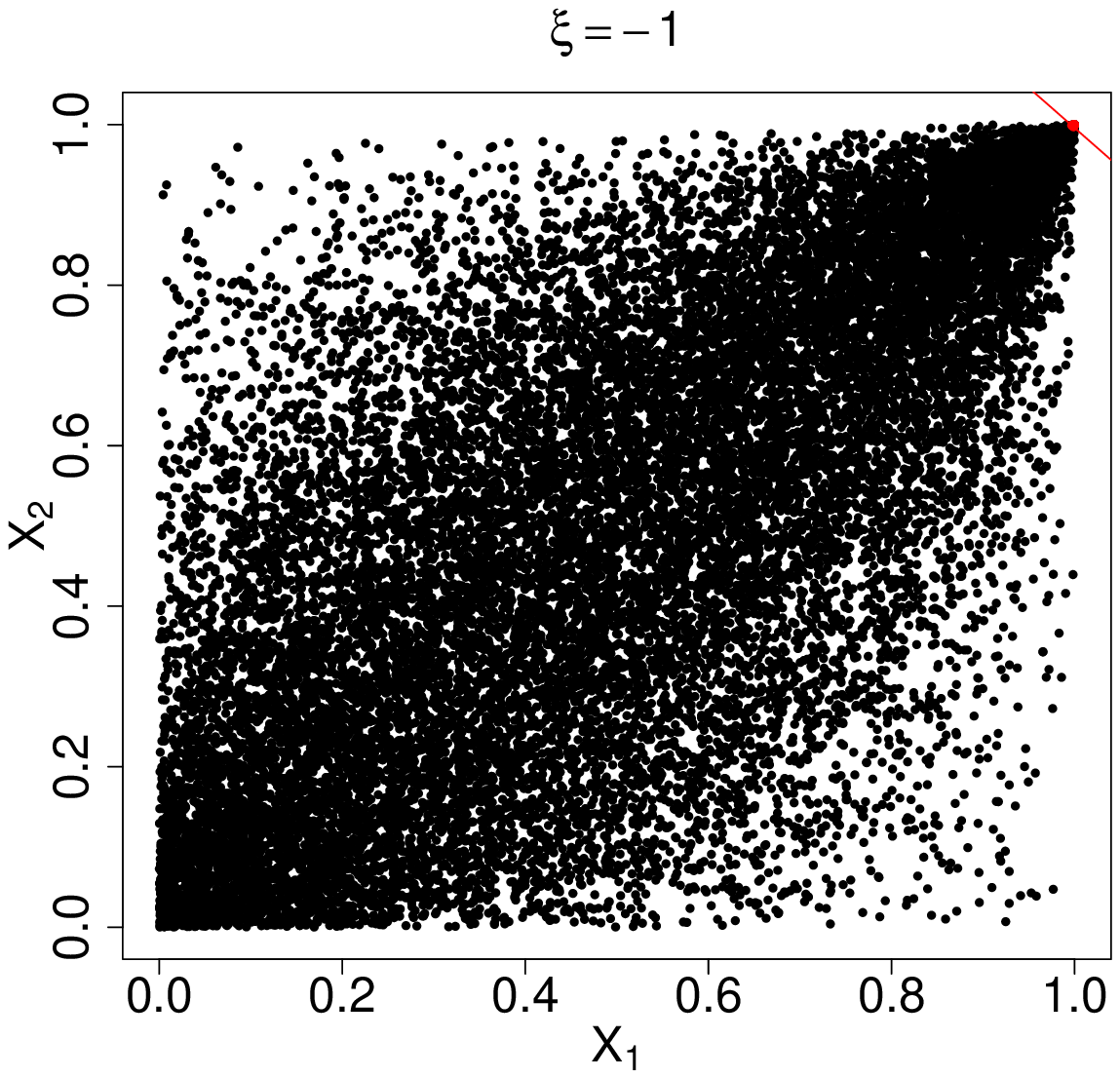} 
\end{minipage}
\begin{minipage}{0.32\linewidth} 
\centering
\includegraphics[width=\linewidth,angle=-90]{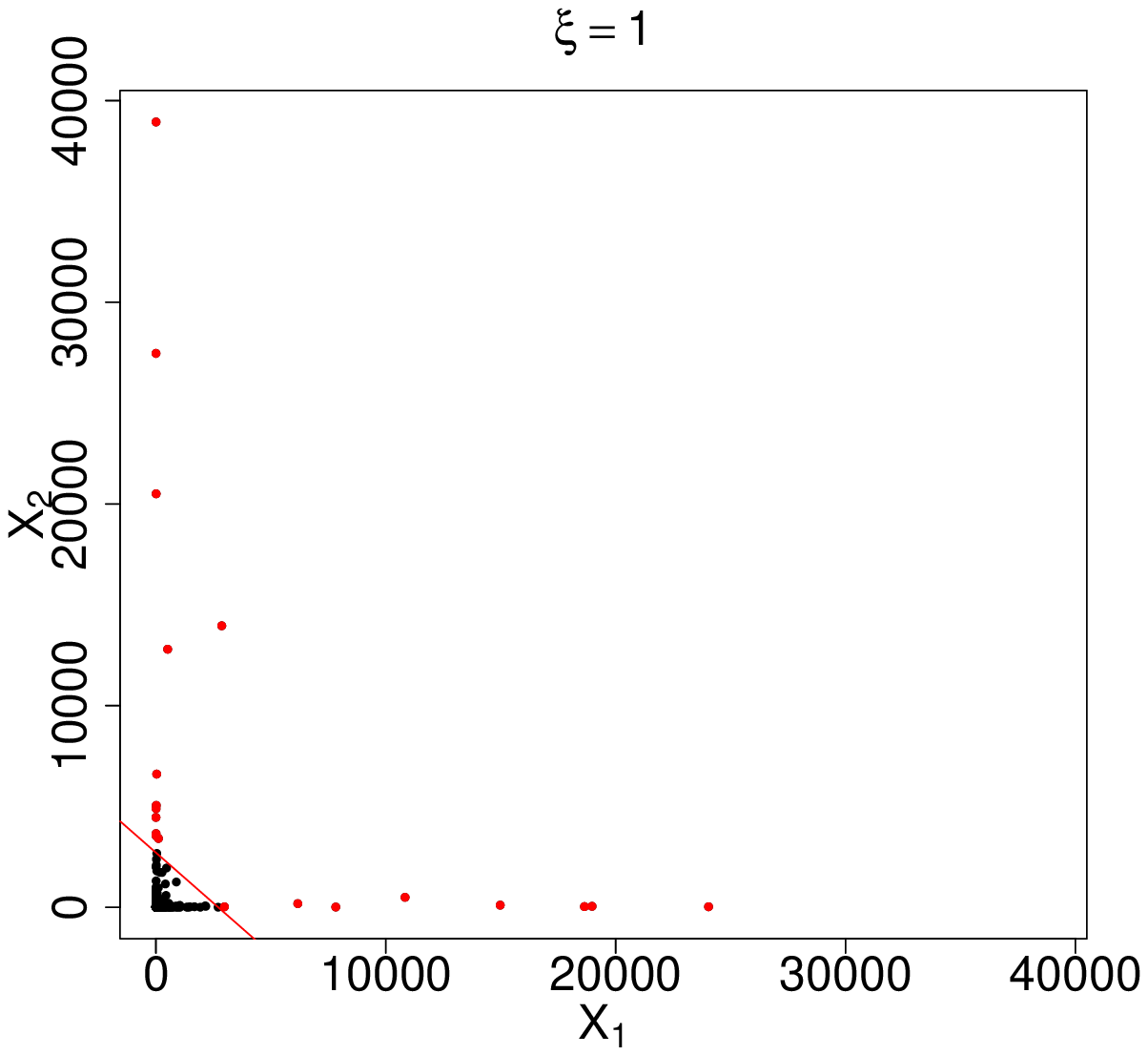} 
\end{minipage}
\begin{minipage}{0.32\linewidth} 
\centering
\includegraphics[width=\linewidth,angle=-90]{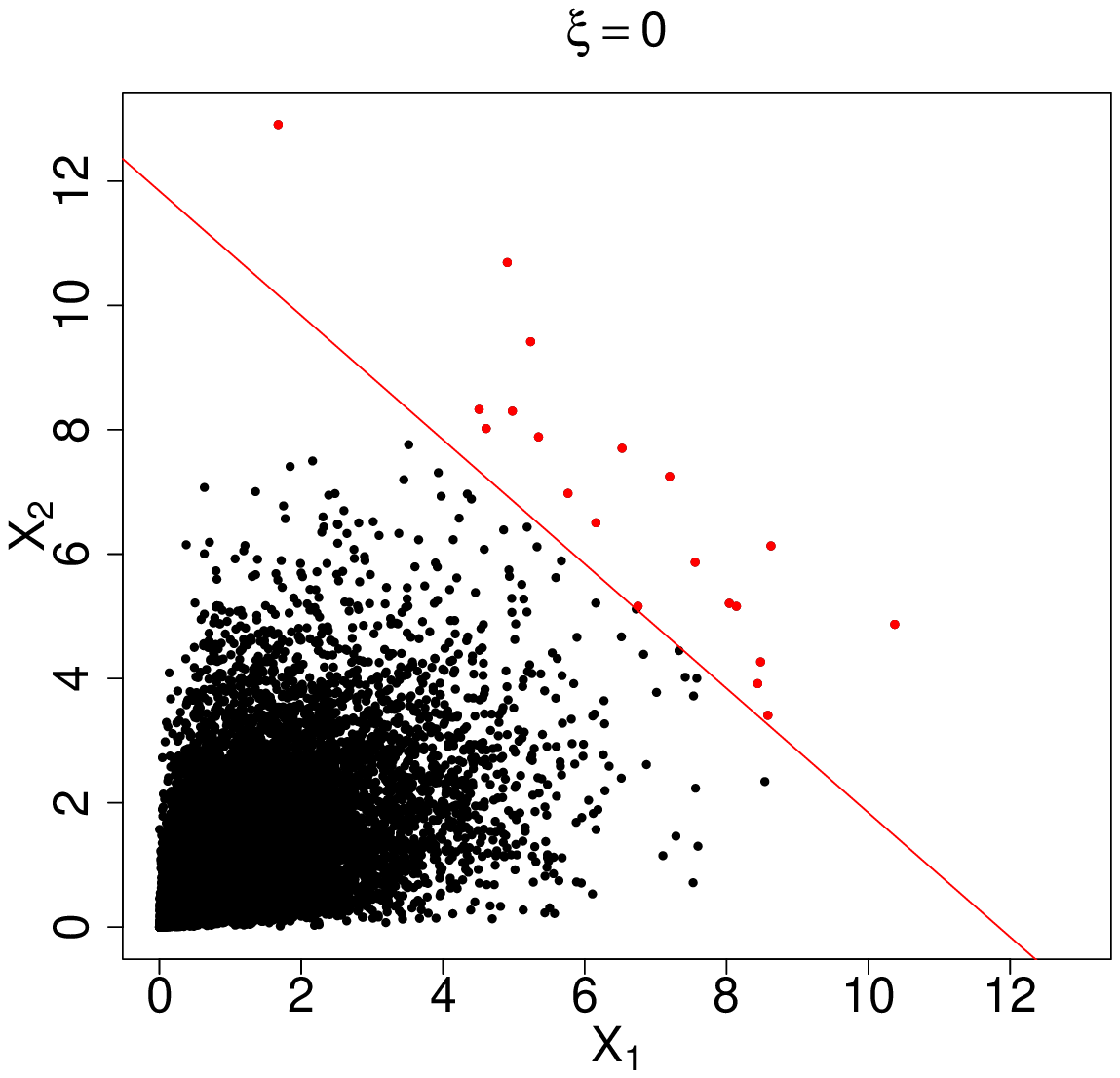}  
\end{minipage}
\begin{minipage}{0.32\linewidth} 
\centering
\includegraphics[width=\linewidth,angle=-90]{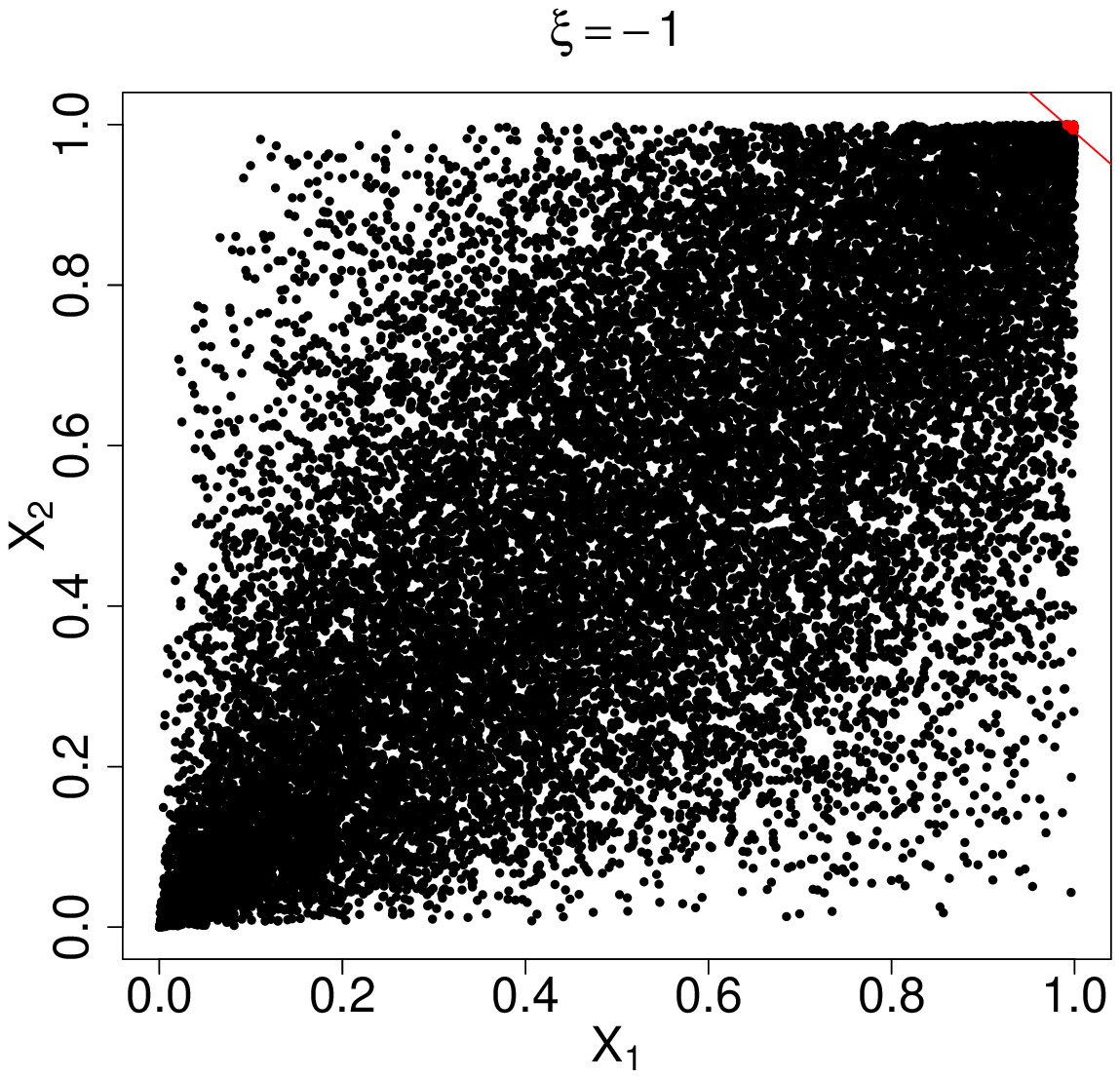} 
\end{minipage}
\caption{Scatter plots of $20000$ simulated $X_1,X_2 \sim \mbox{GPD}(1,\xi)$ with copula \eqref{logistic} (top) and \eqref{invlogistic} (bottom). Both copulae take parameter value $\gamma = 0.5$ and so $(\chi,\eta)=(2-2^{1/2},1)$ and $(0,2^{-1/2})$ in the two rows, respectively. The red points are those for which  $X_1 + X_2>r_{0.999}$, the estimated $0.999-$quantile of $R=X_1+X_2$.}
\label{SimBiv}
\end{figure}
Fig. \ref{SimBiv} motivates our choice of the regions on which we focus for characterising dependence within $(X_1,X_2)$ to derive the extremes of $R$. Here we plot simulated $X_1,X_2 \sim \mbox{GPD}(1,\xi)$ with dependence induced through the logistic and inverted-logistic copula, \eqref{logistic} and \eqref{invlogistic}, respectively. The regions of $(X_1,X_2)$ for which $R \geq r_{0.999}$ are shown, with points in these regions highlighted in red. The combinations of $(X_1,X_2)|(R>r_{0.999})$ highlight which aspects of the copula are important for studying the tail behaviour of $R$. These combinations are similar for different copulae, or dependence structures, but differ for different signs on the marginal shape parameter. For $\xi \leq 0$, the large values of $R$ occur for values which are large in both marginals, which suggests that the important regions of the copula are those where both arguments are simultaneously large; \cite{Ledford1997,Ledford1998} detail dependence in these regions.
Conversely, Fig. \ref{SimBiv} illustrates that for $\xi>0$, large values of $R$ occur when $(X_1,X_2)$ is extreme in at least one component. We thus require a model that considers the distribution of one variable whilst the other is already extreme; which is covered by the characterisation of \cite{Heffernan2004}. We use both approaches for describing limiting dependence of $(X_1,X_2)$ and detail these characterisations in Section~\ref{Setup-sec}.\par
In Section~\ref{intro}, we specified that throughout we would assume
that $X_i>0$ with $X_i \sim \mbox{GPD}(\sigma_i, \xi_i)$ for $i = 1,2$.
These assumptions are clearly highly restrictive when describing marginal behaviour, but as our interest lies in the upper-tail behaviour of $R$, we find that the full distribution of $X_i$ is not always relevant. For example, Fig.~\ref{SimBiv} indicates that when $\max\{\xi_1,\xi_2\} <0$, the combinations of $(X_1,X_2)$ which give large $R$ require both $X_i$ variables to be in their upper-tails. When $\xi_1=\xi_2=\xi  \geq 0$ and $(X_1,X_2)$ are positively dependent in their extremes, large values of $R$ tend to occur when both marginal variables are in their tails. In the case where extremal dependence is weak and the marginal tails are heavy, then $R$ is dominated by only one large marginal variable; the distribution of the values in the body of the smaller variable is not important for the characteristics of the upper-tail of $R$.\par
These arguments indicate that it is predominantly the upper-tail of the marginal variables that are important. The widely adopted approximation for the upper-tails of arbitrary marginal variables is that, for some high quantile $u_i > 0$ of $X_i$, that $(X_i - u_i)|(X_i > u_i)$ follows a GPD \citep{Pickands1975}. Our approach is consistent with this, following the threshold stability property \citep{coles2001} of the GPD: that for all $0<u_i<x_i^F$ we have $(X_i - u_i)|(X_i > u_i)\sim \mbox{GPD}(\sigma_i-\xi_i u_i,\xi_i)$, and so our approach is consistent with the usual tail model without any loss of generality. Thus, our modelling of the marginal distribution has the following properties: it avoids the arbitrary choice of $u_i$; it determines the shape parameter of the tail of $R$ for all $\xi_i$; when $\xi_i=0$ it uniquely determines the scale parameter of the tail through $\sigma_i$; and if the marginal variables are not lower bounded by zero, then similar results are obtained by location shifting the $X_i$, where $X_i$ has a finite lower bound.
 \section{Limit results}
  \label{LimRes-sec}
\subsection{Background and model set-up}
  \label{Setup-sec} 
 In Section~\ref{Results-sec}, we present our results for $\Pr\{R \geq r\}$ in the form
\begin{equation}
\label{Rforms}
\Pr\{R \geq r\} \sim \begin{cases}
K_{1}r^{-1/\xi_R},\;\;&\text{if}\;\; \xi_R > 0,\\
K_{2}\exp\left\{-r/\sigma_R\right\},\;\;&\text{if}\;\; \xi_R = 0,\\
K_{3}\left\{1-r/r^F\right\}^{-1/\xi_R},\;\;&\text{if}\;\; \xi_R < 0,
\end{cases}
\end{equation}
as $r$ tends to $r^F$, the upper-endpoint of $R$, which is infinite if $\xi_R \geq 0$ and is finite when $\xi_R < 0$. Here $\sigma_R > 0$ and $K_{1},K_2,K_{3}>0$ are constants; note that \eqref{Rforms} can be extended by replacing $K_1$ by a slowly varying function, $K_2$ by a regularly varying function, and $K_3$ by a function the converges to a non-zero finite constant as $r\rightarrow r^F$. From expression \eqref{Rforms} it can be seen that the tail of $R$ is predominantly determined by $\xi_R$, with $\sigma_R$ important when $\xi_R = 0$, and $r^F$ when $\xi_R < 0$. Note that in general $r^F \leq x_1^F + x_2^F$, where $x_i^F$ is the upper-endpoint of $X_i$ for $i\in\{1,2\}$, but for the copulae considered in this section the equality holds. In \ref{Rforms-link-sec}, we show how expression \eqref{Rforms} links to the GPD tail formulation which is typically required for modelling using \eqref{GPDCDF}. \par
We now describe the extremal dependence characteristics that we assume for $(X_1,X_2)$. \cite{Ledford1997} present an extension of \eqref{etaEq1} which was extended by \cite{Ramos2009}. Presented here for general marginals $F_1$ and $F_2$, they characterise the joint survival function as
\begin{equation}
\label{landt}
\Pr\left\{F_1(X_1) > 1-1/x_1, F_2(X_2) > 1-1/x_2\right\} \sim \mathcal{L}(x_1+x_2)(x_1x_2)^{-\frac{1}{2\eta}}g\left(x_1/(x_1+x_2)\right),
\end{equation}
for any $x_1 \rightarrow \infty,\; x_2 \rightarrow \infty$ such that $x_1/(x_1+x_2)\rightarrow w$ for $0 < w < 1$, and where $\mathcal{L}(\cdot)$ is a positive slowly-varying function and $g:(0,1)\rightarrow \mathbb{R}_+$ is a continuous function. \cite{Ledford1997} have different powers of $x_1$ and $x_2$ to \eqref{landt} which then requires that $g$ satisfies a property they term quasi-symmetry; however, \cite{Ramos2009} use equal powers of $x_1$ and $x_2$ in the term $(x_1x_2)^{-1/(2\eta)}$ which removes the need for this property. \cite{Ledford1997} provide examples of $g$ for certain copulae, e.g., for the bivariate extreme value logistic copula, they show that $g(w)=\{w(1-w)\}^{-1/2}[1-(w^{1/\gamma}+(1-w)^{1/\gamma})^\gamma]$ for $\gamma$ defined in \eqref{logistic}, and for the inverted logistic copula defined in \eqref{invlogistic}, they show that $g(w)\rightarrow 1$ for all $w\in(0,1)$ as $x_1+x_2\rightarrow \infty$.\par
 \citet{Heffernan2004} and \cite{keef2013} quantify extremal dependence between variables by conditioning on one variable being extreme. We focus on their non-negative association form only. To model extremal dependence in $(X_1,X_2)$, they consider the transformed variables $Y_1=-\log\{1-F_1(X_1)\}$ and $Y_2=-\log\{1-F_2(X_2)\}$, such that $Y_1,Y_2$ are standard exponential random variables. Under the assumption that there exist normalising functions $a:\mathbb{R}\rightarrow \mathbb{R},b:\mathbb{R}\rightarrow \mathbb{R}_{+}$, then for any fixed $z \in \mathbb{R}, y \in \mathbb{R}_+$ and for any sequence $u \rightarrow \infty$, we have
\begin{equation}
\label{heffeq1}
\Pr\left\{[Y_2-a(Y_1)]/{b(Y_1)} < z, Y_1 - u > y | Y_1 > u\right\}\rightarrow \exp(-y)G_Z(z),\;\;\;\text{as}\;\; u \rightarrow \infty,
\end{equation}
 where $G_Z(\cdot)$ is non-degenerate and $\lim_{z\rightarrow \infty}G_Z(z)=1$. Often the normalising functions are simplified to location and scale parameters, i.e., $a(y)=\alpha y$ for $\alpha \in [0,1]$ and $b(y)=y^\beta$ for $0  \leq \beta < 1$. The values of $\alpha$ and $\beta$ determine the strength of dependence between $Y_1$ and $Y_2$, and, thus, between $X_1$ and $X_2$. For example, asymptotic dependence between the two is implied by values $\alpha = 1, \beta = 0$ with $\chi=\int^\infty_0[1-\bar{G}_Z(-z)]\exp(-z)\mathrm{d}z$ for $\bar{G}_Z(z)=1-G_Z(z)$. Within the class of asymptotic independence, we have $\alpha < 1, \beta \geq 0$, with $\alpha = \beta = 0$ giving near perfect independence; we further require $G_Z(\cdot)$ to be standard exponential if $(X_1,X_2)$ are independent.\par
We motivated the use of \eqref{landt} and \eqref{heffeq1} as dependence models in Section~\ref{motivation-sec} by illustrating, through a numerical study, that the largest values of $R$ typically occur when both $X_1$ and $X_2$ are large if $\max\{\xi_1,\xi_2\}\leq 0$, and occur when only one of $X_1,X_2$ is large if $\max\{\xi_1,\xi_2\}> 0$.
\subsection{Results}
  \label{Results-sec}
 We now present the results for the tail behaviour of $R=X_1+X_2$ derived by using the limiting structures described in \eqref{landt} and \eqref{heffeq1} to model dependence in $(X_1,X_2)$.  Recall in \eqref{sum} we define $R$ as a weighted sum, i.e., $R = \omega_1X_1+\omega_2X_2$ with $0 < \omega_1,\omega_2 < 1$ and $\omega_1+\omega_2=1$. By setting $X^*_i=\omega_iX_i$ where $X_i \sim \mbox{GPD}(\sigma_i,\xi_i)$ it follows from \eqref{GPDCDF} that $X^*_i \sim \mbox{GPD}(\omega_i\sigma_i,\xi_i)$ and $R=X^*_1+X^*_2$, and so we present results for $R=X_1+X_2$ without loss of generality. We begin with Theorems~\ref{prophomo} and \ref{propzero}, which detail the cases where the marginal shape parameters are equal and non-zero, and zero, respectively. Theorems~\ref{prophetero} and \ref{propdiffsigns} provide results for the cases where the marginal shapes are unequal; Theorem~\ref{prophetero} covers those cases where both shapes are strictly negative and the other cases are covered by Theorem~\ref{propdiffsigns}.\par
{Assumptions for Theorems 1-3.} We make the assumption that $X_1 \sim \mbox{GPD}(\sigma_1,\xi_1)$ and ${X_2\sim \mbox{GPD}(\sigma_2,\xi_2)}$, with distribution functions defined in \eqref{GPDCDF}, and that the extremal dependence in $(X_1,X_2)$ satisfies the regularity conditions for model \eqref{landt} with $\eta \geq 1/2$ ($\bar{\chi} \geq 0$); we further assume that there exists a fixed $v > 0$ such that, for all $y>v$, $\mathcal{L}(y)$ is a positive constant which is absorbed by the function $g$. We assume that model \eqref{landt} holds in equality for $x_1+x_2\geq \max\{c,u^*\}$ for some fixed constant $c$, $0<c < r^F$, and where $u^*=\max\{x_1^F,x_2^F\}=\max\{-\sigma_1/\xi_1,-\sigma_2/\xi_2\}$ if $\max\{\xi_1,\xi_2\}<0$, and is 0 otherwise; that is, we require that model \eqref{landt} holds for large $R$. We assume that the first- and second-order derivatives of $g$ exists; further assumptions on $g$ are made for specific cases. If $\min\{\xi_1,\xi_2\}=\xi>0$, we require an additional assumption that the limit in \eqref{heffeq1} holds in equality for some fixed $u > 0$ and that the residual distribution $G_Z$ is differentiable. For the theorems that require different conditions for $g$ in \eqref{landt}:\par
{Condition 1.}
There exists a fixed $v^*>0$, such that for $r=x_1+x_2>v^*$, we assume that $g(\omega_x) = 1$ for all $\omega_x=\exp(x_1/\sigma_1)/[\exp(x_1/\sigma_1)+\exp(x_2/\sigma_2)]\in[0,1]$; or equivalently, the density of $X_1$ and $X_2$ factorises when $R=X_1+X_2 >v^*$, and where $X_i \sim \mbox{Exp}(1/\sigma_i)$ for $\sigma_i > 0$ and $i\in\{1,2\}$.\par
{Condition 2.}
The tails of $g$ satisfy $g(w) \sim K_gw^{\kappa}$ as $w\rightarrow 0$ and $g(w) \sim K_g(1-w)^{\kappa}$ as $w\rightarrow 1$ for constant $K_g>0$ and fixed $0 \leq \kappa < 1/(2\eta)$.\par
{Condition 3.}
With $\eta=1$ and as $w \rightarrow 0$ or $w\rightarrow 1$, we have that $g(w) \sim w^{-1/2}(1-w)^{-1/2}[1-H((1-w)^{-1},w^{-1})]$, where the bivariate function $H$ is homogeneous of order $-1$. Moreover, the first, $H_1$ and $H_2$, and second-order, $H_{12}$, partial derivatives of $H$ exist and are continuous, and $\lim_{z\rightarrow\infty}H(z,t)=\lim_{z\rightarrow\infty} H(t,z)=t^{-1}$ for $t>0$. We present two sub-conditions: Condition 3a, $H_{12}(1,z)\sim -K_{H_1}z^{c_1}$ as $z\rightarrow0$ for constants $K_{H_1}>0$ and $c_1>-1$ and $H_{12}(1,z)\sim -K_{H_2}z^{c_2}$ as $z \rightarrow \infty$ for constants $K_{H_2}>0$ and $c_2<-2$; Condition 3b, if $\xi_1 > \xi_2$, then $\lim_{z\rightarrow\infty}H_1(1,z)< \infty$ and $\lim_{z\rightarrow 0}H_1(1,z)= 0$, and otherwise if $\xi_1 < \xi_2$, then $\lim_{z\rightarrow\infty}H_2(z,1)< \infty$ and $\lim_{z\rightarrow\infty} H_2(z,1)=0$.\par
Note that Conditions 2 and 3 are mutually exclusive.
  \begin{theorem}
  \label{prophomo}
 Under the assumptions stated above, if $\xi_1=\xi_2=\xi \neq 0$, then
  \[
  \Pr\{R \geq r\} \sim \begin{cases}
  K\left(1+\frac{\xi r}{\sigma_1+\sigma_2}\right)^{-\frac{1}{\eta\xi}},\;\;\;&\text{if}\;\;\xi < 0,\\
 
   K_u^{*}r^{-1/\xi},\;\;\;&\text{if}\;\;\xi > 0, \;\text{Condition 2 holds or Condition 3 holds},\\
   \end{cases}
  \]
  as $r\rightarrow r^F$, where $r^F=\infty$ for $\xi > 0$ and $r^F = -(\sigma_1+\sigma_2)/\xi$ for $\xi < 0$, and for constants $K$ and $K_u^{*}$ defined in \eqref{negK} and \eqref{posKall}, respectively.
  \end{theorem}
\begin{theorem}
  \label{propzero}
  Under the assumptions stated above, if $\xi_1=\xi_2=0$, then
   \[
   \Pr\{R \geq r\} \sim \begin{cases}
   \frac{\sigma_{max}}{\sigma_{max}-\sigma_{min}}\exp\left(-\frac{r}{2\eta\sigma_{max}}\right),\;\;\;&\text{if}\;\;\sigma_1 \neq \sigma_2, \text{Condition 1 holds},\\
    \frac{r}{2\eta\sigma}\exp\left(-\frac{r}{2\eta\sigma}\right),\;\;\;&\text{if}\;\sigma_1 = \sigma_2=\sigma, \text{Condition 1 holds},\\
      K\exp\left(-\frac{r}{(\sigma_1+\sigma_2)}\right),\;\;\;&\text{if}\;\; \text{Condition 3a holds,}\\
      \end{cases}
      \]
      as $r\rightarrow \infty$ and for constant $K$ defined \eqref{Isimeq}, where $\sigma_{max}=\max\{\sigma_1,\sigma_2\}$ and $\sigma_{min}=\min\{\sigma_1,\sigma_2\}$.
  \end{theorem}
\hspace{-0.5cm}Note that there is a power term in the second case for $\Pr\{R \geq r\}$ given by Theorem~\ref{propzero} that is not covered by the general form given by \eqref{Rforms}. However, the form in Theorem~\ref{propzero} is in the domain of attraction of a GPD with shape and scale parameters zero and $2\eta\sigma$, respectively.
  \begin{theorem}
  \label{prophetero}
   Under the assumptions stated above, if $\xi_1 \neq \xi_2$ and $\max\{\xi_1,\xi_2\} < 0$, then
  \[
  \Pr\{R \geq r\} \sim \begin{cases}
  K_{(2)}\left(1+\frac{\xi_1\xi_2 r}{\sigma_1\xi_2+\sigma_2\xi_1}\right)^{-\frac{1}{\xi_{max}}\left(\frac{1}{2\eta}+\kappa\right)-\frac{1}{\xi_{min}}\left(\frac{1}{2\eta}-\kappa\right)},\;\;\;&\text{if Condition 2 holds},\\
  K_{(3b)}\left(1+\frac{\xi_1\xi_2 r}{\sigma_1\xi_2+\sigma_2\xi_1}\right)^{-\frac{1}{\xi_{max}}},\;\;\;&\text{if Condition 3b holds},\\
   \end{cases}
  \]
  as $r\rightarrow r^F = -(\sigma_1/\xi_1+\sigma_2/\xi_2)$, and for constants $K_{(2)}>0$ and $K_{(3b)}>0$ defined in \eqref{heteronegK} and \eqref{heteronegK2}, respectively, and $\xi_{max}=\max\{\xi_1,\xi_2\},\; \xi_{min}=\min\{\xi_1,\xi_2\}$. 
  \end{theorem}
\hspace{-0.5cm}The set conditions on the dependence in $(X_1,X_2)$ given for Theorems~\ref{prophomo}-\ref{prophetero} are not necessary for Theorem~\ref{propdiffsigns}. Whilst Theorems~\ref{prophomo}-\ref{prophetero} only apply when $\bar{\chi}\geq 0$, i.e., non-negative association in $(X_1,X_2)$, Theorem~\ref{propdiffsigns} applies for any $\bar{\chi} \in (-1,1]$.
   \begin{theorem}
  \label{propdiffsigns}
 If $\xi_1 \neq \xi_2$ and $\max\{\xi_1,\xi_2\} \geq 0$, then
  \[
  \Pr\{R \geq r\} \sim \begin{cases}
 \left(\xi_{max}/\sigma_{max}\right)^{-1/\xi_{max}}r^{-1/\xi_{max}},\;\;\;&\text{if}\;\;\max\{\xi_1,\xi_2\} >0,\\
   C\exp\left(-r/\sigma_{max}\right),\;\;\;&\text{if}\;\;\max\{\xi_1,\xi_2\} = 0,\\
   \end{cases}
  \]
 as $r\rightarrow \infty$ and where $\xi_{max}=\max\{\xi_1,\xi_2\}$ and $\sigma_{max}=\{\sigma_i; i\text{ is s.t. } \xi_i = \xi_{max}\}$ and for constant $C\in[1,C_D^*]$ with $C_D^*$ defined in \eqref{C1eq}.
  \end{theorem}
Using a different approach \cite{koutsoyiannis2020stochastics} provides a similar result to Theorem~\ref{propdiffsigns} when $\min\{\xi_1,\xi_2\}>0$. We further note that the cases where $\min\{\xi_1,\xi_2\}>0$ in Theorems~\ref{prophomo} and \ref{prophetero} agree with Breiman's Lemma \citep{breiman1965some}, as we have $\xi_R = \max\{\xi_1,\xi_2\}$.
  \section{Copula examples}
  \label{cop-sec}
We now compare the limit results detailed in Section~\ref{Results-sec} with results for the upper-tail behaviour of $R$ when dependence in $(X_1,X_2)$ is fully modelled using copula families and their marginal models remain the same, i.e., $X_i \sim \mbox{GPD}(\sigma_i,\xi_i)$ for $i\in\{1,2\}$. The assumptions we made in Section~\ref{Results-sec} hold in some cases and in these we obtain identical results to Section~\ref{Results-sec}. However, where the assumptions of Section~\ref{Results-sec} are too strong, our direct derivations from the copulae, with details in \cite{richardsthesis}, provide insight into the tails of $R$ in these specific cases and they show some features of Theorems~\ref{prophomo}-\ref{prophetero} still hold even when their conditions fail to hold. We consider the extreme value copula and the inverted extreme value copula and the limiting forms of these two classes, i.e., perfect dependence and independence. We further consider a standard Gaussian copula with correlation parameter $\rho\in(0,1)$.\par
  The extreme value copula takes the form \begin{equation}
  \label{bev}
  C_{ev}(u,v)=\exp \left\{-V(-1/\log(u),-1/\log(v))\right\},
  \end{equation}
  where $V(x,y)=2\int^1_0\max\left\{w/x,(1-w)/y\right\}\mathrm{d}M(w)$ is a homogeneous function of order $-1$ and $M(w)$ is a univariate distribution function/probability measure for $w \in [0,1]$, which has expectation $1/2$. Note that $1 \leq V(1,1) \leq 2$, where the boundary cases correspond to special cases of the extreme value copula, i.e., we have perfect dependence, and independence, between $X_1$ and $X_2$ when $V(1,1)=1$ and $V(1,1)=2$ respectively. This copula gives $\eta = 1$ ($\chi=2-V(1,1), \;\bar{\chi}=1$) and $\eta = 1/2$ ($\chi=0, \;\bar{\chi}=0$) when $V(1,1) <2$ and $V(1,1)=2$, respectively; that is, this copula exhibits asymptotic dependence, or independence, only. Furthermore, \cite{Ledford1997} illustrate that this copula satisfies Condition 3a/3b, with $H = V$ and $\kappa = 1/2$, needed for Theorems~\ref{propzero} and \ref{prophetero}. \par
  The inverted extreme value copula follows by inverting \eqref{bev}, see \cite{Ledford1997}, and is defined through its survival copula 
  \begin{equation}
  \label{Ibev}
  \bar{C}_{iev}\{u,v\}=\exp \left\{-V(-1/\log(1-u),-1/\log(1-v))\right\},
  \end{equation} 
with a similarly defined $V$. This, and the Gaussian copula, have $\eta =V(1,1)^{-1}$ and $\eta=(1+\rho)/2$, respectively, where $\bar{\chi}=2\eta - 1$ and $\chi = 0$ for both copulae. The bivariate extreme value logistic \eqref{logistic}, and inverted logistic \eqref{invlogistic}, copulae are subclasses of \eqref{bev} and \eqref{Ibev}, respectively.
When discussing results pertaining to copulae \eqref{bev} and \eqref{Ibev}, we assume that the first- and second-order partial derivatives of $V$ exist, i.e., the existence of a joint density, hence this excludes perfect dependence which is derived separately.  \par
We report the parameters that determine the leading order behaviour of $\Pr\{R \geq r\}$ as $r \rightarrow \infty$ as given by form \eqref{Rforms}, i.e., $\xi_R \neq 0$ and when $\xi_R=0$ we report $\sigma_R$. We partition the space of all $(\xi_1,\xi_2)$ into three cases: (i)~$\xi_1 = \xi_2 = 0$, (ii)~$\max\{\xi_1,\xi_2\} <0$ with $\xi_1 \neq \xi_2$ and (iii)~all other possibilities not covered by cases (i) and (ii). Case (iii) covers a large number of sub-cases, e.g., all those covered by Theorems~\ref{prophomo} and \ref{propdiffsigns}, and no further insight into the upper-tails of $R$ is revealed when modelling dependence in $\mathbf{X}$ using copulae; \cite{richardsthesis} finds the same results as detailed in Section~\ref{Results-sec} suggesting that, for this case, modelling dependence using the limiting models of \cite{Ledford1997} and \cite{Heffernan2004} is sufficient to derive the first-order behaviour of the upper-tail of $R$.
\par
However, for cases (i) and (ii), then \cite{richardsthesis} shows that some further insight into the upper-tail behaviour of $R$, relative to that provided by Theorems~\ref{propzero} and \ref{prophetero}, can be gained by modelling dependence with copulae. In these cases, we report the differences between the results derived by \cite{richardsthesis} and those covered by Theorems~\ref{propzero} and \ref{prophetero} in Table~\ref{cop-pars}. Observe that for the classical scenarios of asymptotic dependence $(\chi >0)$ and independence $(\bar{\chi}=0)$, then both sets of results agree; the parameters $\xi_R$ and $\sigma_R$ are the same whether derived from first principles using the copula or obtained by Theorems~\ref{prophomo} and \ref{propzero} under conditions that do not fully apply to the specific copula. However, for asymptotic independence, i.e., $0 \leq \bar{\chi} < 1$, then the results only agree in the limit as $\bar{\chi} \downarrow 0$. As illustrated by \cite{Ledford1997}, none of Conditions 1-3 are met by either the inverted extreme value, or the Gaussian, copulae. For these copulae, the parameters in \eqref{Rforms} cannot be represented as the product of a function of the marginal parameters and the summary measure $\eta$; instead, the upper-tail behaviour of $R$ is driven by a function of both the marginal parameters and dependence structure which cannot be factorised, indicating a more subtle relationship between the marginal shapes, extremal dependence structure and upper-tail behaviour of $R$.
\begin{sidewaystable}
 \centering
 \caption{Parameter values for $R=X_1+X_2$ where $(X_1,X_2)$ have GPD margins with either $\xi_1=\xi_2=0$ or $\xi_1\neq\xi_2$ with $\max\{\xi_1,\xi_2\} < 0$, and for $h(w)=\sigma_1w -2\rho\sqrt{\sigma_1\sigma_2w(1-w)}+\sigma_2(1-w)$, for $0 < w <1, \min\{\sigma_1,\sigma_2\} > 0$ and $0 \leq \rho < 1$. All other cases for the GPD margins are discussed in text in Section~\ref{cop-sec}.}
 \label{cop-pars}
 \hrule
 \begin{tabular}{lcccc} 
&&$\xi_1=\xi_2=0\Rightarrow\xi_R=0$&&$\max\{\xi_1,\xi_2\} <0,\; \xi_1 \neq \xi_2$\\
\cline{3-3}\cline{5-5}
Dependence Structure &&$\sigma_R$&&$\xi_R$\\
 \cline{1-1}\cline{3-3}\cline{5-5}
  Theorems~\ref{propzero}/\ref{prophetero}, $\chi > 0$& & $\sigma_1+\sigma_2$&& $\max\{\xi_1,\xi_2\}$\\
    Theorems~\ref{propzero}/\ref{prophetero}, $0\leq\bar{\chi} < 1$ && $2\eta\max\{\sigma_1,\sigma_2\}$&&$2\eta\left(\xi_1^{-1}+\xi_2^{-1}\right)^{-1}$\\
Independence, $\bar{\chi} = 0$ &&$\max\{\sigma_1,\sigma_2\}$&&$\left(\xi_1^{-1}+\xi_2^{-1}\right)^{-1}$\\
Perfect dependence, $\chi = 1$&&$\sigma_1+\sigma_2$&& $\max\{\xi_1,\xi_2\}$\\
Extreme value copula, $\chi> 0$&&$\sigma_1+\sigma_2$&& $\max\{\xi_1,\xi_2\}$\\
\shortstack{Inverted extreme\\value copula, $0\leq\bar{\chi} < 1$} &&$\max\limits_{0 \leq w \leq 1}\{\left[V\left(\sigma_1/w,\sigma_2/(1-w)\right)\right]^{-1}\}$&&$-1/V(-\xi_1,-\xi_2)$\\
\shortstack{Standard Gaussian,\\ $\bar{\chi} = \rho,\rho\in[0,1)$}&&$(1-\rho^2)\max\limits_{0 \leq w \leq 1}\{h(w)^{-1}\}$&&$(1-\rho^2)\{\xi_1^{-1}+2\rho(\xi_1\xi_2)^{-1/2}+\xi_2^{-1}\}^{-1}$
 \end{tabular}
 \hrule
 
 \end{sidewaystable}
 \section{Application to aggregated environmental data}
\label{Sec-Appl}
\subsection{Introduction}
We now present an application of the results discussed in Section \ref{Results-sec} to climate model data. We study precipitation and temperature data, which have heavy and bounded marginal upper-tails respectively. Both datasets are obtained from the UK climate projections 2018 (UKCP18) \citep{ukcp18report} which contains values aggregated over a given time interval and a spatial grid-box. The size of these grid-boxes and the specified time interval differ between the two studies. In both cases, we investigate the marginal upper-tail for the variables observed at a configuration of grid-boxes and the spatial average of them over adjacent boxes. \par
Recall from Section \ref{Setup-sec} that the driving factor for the extremal behaviour of the aggregates is the GPD shape parameter, $\xi$. We focus on just the relationship between estimates of $\xi$ for the marginal variables and $\xi$ for the aggregates. To investigate this relationship, we begin with a $2$ by $2$ configuration of adjacent grid-boxes. For each pair of adjacent grid-boxes, we conduct inference on three variables. Marginally, for each grid-box, we fit the GPD to excesses above the sample $p$-th quantile using maximum likelihood methods, under the assumption that observations are independent and identically distributed \citep{coles2001}. Following many spatial extreme value applications \citep{Coles1990,Coles1996,Fowler2003,Coelho2008,Li2019, Davison2012,davison2019}, we anticipate that the shape parameters for adjacent grid-boxes should be identical. Therefore we also pool information across grid-boxes with a model that the distribution of excesses in adjacent grid-boxes $i$ and $j$ are GPD$(\sigma_i,\xi)$ and GPD$(\sigma_j,\xi)$, respectively, i.e., a common shape parameter but with the scale parameter unconstrained; note that we do not assume a common shape for all four grid-boxes, we make the weaker assumption that only adjacent grid-boxes share a common shape parameter. We then take the spatial aggregate of the data across adjacent grid-boxes at each separate time interval and fit a GPD to excesses of these data above its empirical $p$-th quantile. Quantiles are estimated separately for marginal, pooled and aggregate variables.
 To account for strong spatial and temporal dependence in the data, standard errors for $\xi$ are estimated using a stationary bootstrap \citep{statboot} with $1000$ samples, with temporal block size drawn randomly from a Geometric distribution with expectation corresponding to a week of observations. 
\subsection{Precipitation}
\label{sec-rain}
The data are precipitation flux (mm/day) from a convection permitting model on $2.2\times 2.2 km^2$ grid-boxes and hourly intervals. To account for seasonality, we use only winter, December to February, observations between the years 1980 and 2000. We study a $2 \times 2$ configuration of grid-boxes centred around $(52.18^\circ, 0.14^\circ)$, approximately Cambridge, UK; this is a flat area so no orographic features are important and marginal distributions are expected to be nearly homogeneous. We conduct our analysis on outputs of the model at two spatial resolutions - high, using data on $(2.2)^2km^2$, and coarse, from $(22)^2km^2$ grids. The latter is produced by taking the spatial average over $10$ by $10$ configurations of the former data. We analyse both resolutions to investigate the effect of extremal dependence on the observed results. This is quantified using the measure $\eta$, given in \eqref{landt}, which is estimated as in \cite{coles1999}. All GPD models are fitted to exceedances above $99.5\%$ quantiles. \par
Table \ref{rainstudy} presents estimates and the $95\%$ confidence intervals for the shape parameters using the three inference methods. The marginal shape parameter estimates are predominately positive which suggests that Theorem \ref{prophomo} is relevant, i.e., for a homogeneous marginal shape parameter $\xi>0$, the shape parameter of the aggregate is also $\xi$, regardless of the dependence structure. We aim to see if this applies in the observed tail. 
\begin{table}[h]
\centering
\caption{High resolution precipitation case study: shape parameter estimates and 95\% confidence intervals for margins (black), pooled marginals (red) and aggregate variable (blue).}
 \label{rainstudy}
 \hrule
 \begin{tabular}{lcccc} 
&&\multicolumn{3}{c}{Marginal}\\
\cline{3-5}
&&0.210\;(0.045, 0.339)&&0.197\;(0.037, 0.350)\\
 \cline{3-3}\cline{5-5}
 \multirow{4}{*}{Marginal} &\multirow{2}{*}{0.154\;(-0.030, 0.286)}& \textcolor{red}{0.172\;(0.017, 0.306)}&&\textcolor{red}{0.178\;(0.019, 0.320)}\\
 && \textcolor{blue}{0.160\;(-0.006, 0.288)}&&\textcolor{blue}{0.172\;(0.020, 0.328)}\\
&\multirow{2}{*}{0.225\;(0.040, 0.344)}& \textcolor{red}{0.214\;(0.049, 0.333)}&&\textcolor{red}{0.168\;(-0.001, 0.283)}\\
&& \textcolor{blue}{0.177\;(0.036, 0.316)}&&\textcolor{blue}{0.184\;(0.041, 0.347)}\\
 \end{tabular}
\hrule
\end{table}
Table \ref{rainstudy} shows the point estimates of confidence intervals for $\xi$ using the marginal variables and the pooled analysis. As we observe similar estimates for $\xi$ as well as substantial overlap in the confidence intervals, this suggests that it is reasonable to assume homogeneous marginal shape parameters. Using the same criteria as above, the marginal estimates also have good agreement with $\xi$ for the aggregate variable, suggesting that the positive shape result in Theorem \ref{prophomo} holds well for these data. Pairwise $\eta$ estimates for Table~\ref{rainstudy} fall in the range $[0.956,0.967]$, which suggests strong extremal dependence between the marginal variables.
\begin{table}[h]
\centering
\caption{Coarse resolution precipitation case study: shape parameter estimates and 95\% confidence intervals for margins (black), pooled marginals (red) and aggregate variable (blue).}
 \label{rainstudy2}
 \hrule
 \begin{tabular}{l cccc} 
&&\multicolumn{3}{c}{Marginal}\\
\cline{3-5}
&&0.146\;(-0.033, 0.277)&&0.089\;(-0.024, 0.197)\\
  \cline{3-3}\cline{5-5}

 \multirow{4}{*}{Marginal} &\multirow{2}{*}{0.104\;(-0.083, 0.218)}& \textcolor{red}{0.108\;(-0.015, 0.239)}&&\textcolor{red}{0.101\;(0.000, 0.186)}\\
 && \textcolor{blue}{0.177\;(-0.095, 0.318)}&&\textcolor{blue}{0.011\;(-0.123, 0.085)}\\
&\multirow{2}{*}{0.068\;(-0.055, 0.183)}& \textcolor{red}{0.105\;(-0.119, 0.212)}&&\textcolor{red}{0.082\;(-0.012, 0.176)}\\
&& \textcolor{blue}{0.085\;(-0.082, 0.182)}&&\textcolor{blue}{0.061\;(-0.065, 0.189)}\\

 \end{tabular}
\hrule
\end{table}
To investigate the effect of weaker dependence on the relationship between the marginal and aggregate $\xi$ parameter, we now consider the coarse resolution data and conduct the same analyses as previously; pairwise $\hat{\eta}$ for the coarser data are in the range $[0.859,0.895]$, which is lower than the estimates for Table~\ref{rainstudy}.
Table \ref{rainstudy2} suggests that it is reasonable to assume homogeneous marginal shape parameters at this coarse resolution, as we again observe good agreement between the $\xi$ estimates for both the marginal and pooled variables. We also observe good agreement between $\xi$ for the pooled variables and aggregate variables even with weaker extremal dependence. 
\subsection{Temperature}
The data are average daily temperature ($^\circ$C) from a global climate model scaled to $60\times 60 km^2$ grid-boxes and to account for seasonality we use only summer, June to August, observations. The model is run through the years 1899 to 2099, providing $18000$ observations per grid-box. We consider a $2 \times 2$ configuration of grid-boxes centred around $(53.14^\circ, -1.70^\circ)$, south of the Peak District, UK. As in Section \ref{sec-rain}, we conduct our analyses on outputs of the model at two spatial resolutions - high using data on $(60)^2km^2$ and coarse $(300)^2km^2$ grids; the latter produced by taking the spatial average over $5$ by $5$ configurations of the former data. All GPD models are fitted to exceedances above $98\%$ quantiles.\par
Table \ref{heatstudyshape} presents estimates and the $95\%$ confidence intervals for the shape parameter for the marginal and pooled variables, which suggest that these variables have bounded upper-tails. As such, we consider the results in Theorem \ref{prophomo}; this states that, asymptotically, the shape parameter of the aggregate should be $\eta\xi$ given that the marginal variables have equal, negative shape $\xi<0$. To see if this result is consistent with the observed tails, Table \ref{heatstudyshape} presents estimates and $95\%$ confidence intervals for a scaling of the aggregate shape parameter by $1/\hat{\eta}$, where the estimate $\hat{\eta}$ of $\eta$ is calculated for each bootstrap sample of the aggregate; if Theorem \ref{prophomo} holds for these data, then this should be equal to the marginal $\xi$. 
\begin{table}[h]
\centering
\caption{High resolution temperature case study: shape parameter estimates and 95\% confidence intervals for margins (black) and pooled variable (red). Blue confidence intervals are for a scaling of the aggregated shape parameter by $1/\eta$.}
 \label{heatstudyshape}
 \hrule
 \begin{tabular}{l cccc} 

&&\multicolumn{3}{c}{Marginal}\\
\cline{3-5}
&&-0.156\;(-0.276, -0.067)&&-0.211\;(-0.308, -0.108)\\
   \cline{3-3}\cline{5-5}

 \multirow{4}{*}{Marginal} &\multirow{2}{*}{-0.198\;(-0.310, -0.106)}& \textcolor{red}{-0.180\;(-0.268, -0.106)}&&\textcolor{red}{-0.199\;(-0.293, -0.133)}\\
 && \textcolor{blue}{-0.214\;(-0.339, -0.103)}&&\textcolor{blue}{-0.201\;(-0.318, -0.103)}\\
&\multirow{2}{*}{-0.148\;(-0.266, -0.082)}& \textcolor{red}{-0.165\;(-0.255, -0.069)}&&\textcolor{red}{-0.161\;(-0.250, -0.094)}\\
&& \textcolor{blue}{-0.166\;(-0.278, -0.083)}&&\textcolor{blue}{-0.160\;(-0.297, -0.067)}
 \end{tabular}
 \hrule
\end{table}
Table \ref{heatstudyshape} suggests that we can assume homogeneous marginal shape parameters and these estimates also have clear agreement with the scaled shape parameter for the aggregate variable, suggesting that the negative shape result in Theorem \ref{prophomo} holds well for these data. Pairwise $\eta$ estimates for Table~\ref{heatstudyshape} fall in the range $[0.918,0.981]$, which suggests strong extremal dependence between the marginal variables, and so we repeat the analyses with the coarser data to investigate the effect of weaker dependence on the aggregate shape parameter.
\begin{table}[h]
\centering
\caption{Coarse resolution temperature case study: shape parameter estimates and 95\% confidence intervals for margins (black) and pooled variable (red). Blue confidence intervals are for a scaling of the aggregated shape parameter by $1/\eta$.}
 \label{heatstudyshape2}
 \hrule
 \begin{tabular}{lcccc} 
&&\multicolumn{3}{c}{Marginal}\\
\cline{3-5}
&&-0.113\;(-0.277, -0.020)&&-0.207\;(-0.298, -0.132)\\
  \cline{3-3}\cline{5-5}

 \multirow{4}{*}{Marginal} &\multirow{2}{*}{-0.183\;(-0.280, -0.106)}& \textcolor{red}{-0.145\;(-0.219, -0.088)}&&\textcolor{red}{-0.200\;(-0.272, -0.158)}\\
 && \textcolor{blue}{-0.200\;(-0.356, -0.102)}&&\textcolor{blue}{-0.204\;(-0.342, -0.129)}\\
&\multirow{2}{*}{-0.053\;(-0.317, 0.057)}& \textcolor{red}{-0.138\;(-0.258, -0.067)}&&\textcolor{red}{-0.066\;(-0.255, 0.010)}\\
&& \textcolor{blue}{-0.083\;(-0.338, 0.011)}&&\textcolor{blue}{-0.178\;(-0.410, -0.071)}
 \end{tabular}
 \hrule
\end{table}
Table \ref{heatstudyshape2} suggests that it is still reasonable to assume homogeneous marginal shape parameters at the coarser resolution, as we again observe good agreement between the $\xi$ estimates for both the marginal and pooled variables. We found that pairwise values of $\hat{\eta}$ for Table~\ref{heatstudyshape2} were in the range $[0.789,0.921]$, which suggests weaker extremal dependence than that observed for the high resolution temperature data. We also observe good agreement between these estimates and the estimates for the scaled aggregate shape parameter, confirming that the result in Theorem \ref{prophomo} applies well, even for weaker extremal dependence.
\section{Discussion}
\label{discuss-sec}
In Section \ref{Results-sec}, we provided results that explore the extremal behaviour of $R$; the bivariate aggregate of two GPD random variables, $X_1$ and $X_2$. These results focus primarily on the effect of the marginal shape parameters and dependence within $(X_1,X_2)$ on the shape parameter $\xi_R$ of the aggregate, or its scale parameter if we have $\xi_R = 0$. Through Section \ref{Setup-sec}, we illustrate that the value of the maximum of the marginal shape parameters is generally the most important driver in the tail behaviour of the aggregate.
\par
The results given in Section~\ref{Results-sec} were derived by modelling the dependence in $(X_1,X_2)$ using the limiting extremal dependence models of \cite{Ledford1996} and \cite{Heffernan2004}, whereas results using full copula models are given in Section~\ref{cop-sec}. There is broad agreement between results derived using the two methods, and so we conclude that the extremal behaviour of $R$ is mostly driven by the limiting behaviour of $(X_1,X_2)$ as $x_1 \rightarrow \infty$ and/or $x_2 \rightarrow \infty$, and that modelling the full dependence in $(X_1,X_2)$ is not necessary to capture the first order behaviour of $\Pr\{R \geq r\}$ as $r \rightarrow \infty$.\par
Although we define $R$ for any $d\in\mathbb{N}$ in \eqref{sum}, we constrain the focus of our study to $d=2$. To extend the results in Section~\ref{Results-sec} to $d >2$, we require extensions of characterisations \eqref{landt} and \eqref{heffeq1} to greater dimensions; such variants exist, see \cite{Eastoe2012}, and so it is reasonable to assume that the results given in Section~\ref{Results-sec} can be extended to $d>2$. For example, \cite{richardsthesis} details extensions of some of the copula results given in Section~\ref{cop-sec} to $d>2$. However, we note
that the framework for the proofs in \ref{proof_sec} may not be applicable when deriving results for $d > 2$, as we would require evaluation of $d$-dimensional integrals which may not be feasible in an analytical setting.\par
To illustrate some of the practical utility of the results in the paper, in Section~\ref{Sec-Appl} we undertake inference on the upper-tail behaviour of aggregated precipitation and temperature data, which have heavy and bounded marginal upper-tails, respectively. We aggregate the data as we are interested in the extremal behaviour of the climate processes at lower resolutions; for precipitation, this is for the reasons described in Section \ref{intro}, and for temperature, we are interested in the average extreme heat over a large spatial domain since a heatwave has societal impact owing to it affecting a spatial region not simply a single location. Both datasets are obtained from the UK climate projections 2018 (UKCP18) \citep{ukcp18report} which contains values aggregated over a spatial grid-box; we conduct our analyses using both a fine, and coarse, spatial resolution, as we find weaker extremal dependence exhibited by the latter data. We estimate $\xi_1,\xi_2$ and find that we can reasonably assume $\xi_1=\xi_2$ ($=\xi$ say) as their estimated confidence intervals overlap. We further estimate $\xi_R$ and illustrate that these estimates agree with the results presented in Theorem~\ref{prophomo}; namely that $\xi_R=\xi$ and $\xi_R=\eta\xi$ if $\xi>0$, and $\xi<0$, respectively.
\section*{Acknowledgements} J. Richards gratefully acknowledges funding through the STOR-i Doctoral Training Centre and Engineering and Physical Sciences Research Council (grant EP/L015692/1).  The authors are grateful to the UK Met Office for data and to Simon Brown and Jennifer Wadsworth for supportive discussions.
\appendix
\begin{appendix}
\section{Proofs}
\label{proof_sec}
\textbf{Proof of Theorem \ref{prophomo}.}
  Negative Shape Case ($\xi < 0$). We begin by deriving the joint density of $(X_1,X_2)$ given in \eqref{landt} for GPD margins. We transform $(X_1,X_2)\rightarrow (R,W)$, where $R=X_1 + X_2$ and $W$, an auxiliary variable, and integrate out $W$ to give the density of $R$ and derive its survival function. Combining \eqref{landt} and \eqref{GPDCDF} with $\xi_1=\xi_2=\xi <0$, we have that
  \[
 \Pr\left\{1-\left(1+\xi\frac{X_1}{\sigma_1}\right)^{-1/\xi} > 1-\frac{1}{x_1}, 1-\left(1+\xi\frac{X_2}{\sigma_2}\right)^{-1/\xi} > 1-\frac{1}{x_2}\right\} = \frac{\mathcal{L}(x_1+x_2)}{(x_1x_2)^{\frac{1}{2\eta}}}g\left(\frac{x_1}{x_1+x_2}\right),
  \]
  as $x_1,x_2\rightarrow \infty$ such that the limit of $x_1/(x_1+x_2)$ is bounded by $(0,1)$. Under the assumption that $\mathcal{L}(y)$ acts as a constant which can be absorbed by $g$ for $y > v$ for some $v>0$, we have ${\Pr\left\{X_1>x_1, X_2>x_2\right\} \sim \tilde{x}_1^{-\frac{1}{2\eta\xi}}\tilde{x}_2^{-\frac{1}{2\eta\xi}}g\left(\omega_x\right)}$ for $x_1\rightarrow x_1^F$ and $x_2 \rightarrow x_2^F$, such that $\omega_x= \tilde{x}^{1/\xi}_1/(\tilde{x}^{1/\xi}_1+\tilde{x}^{1/\xi}_2)\rightarrow \omega^*_x \in (0,1)$ and where $\tilde{x}_i=(1+\xi x_i/\sigma_i)$ for $i\in\{1,2\}$. Assuming that the first and second derivatives of $g$ exist, then the density of $(X_1,X_2)$ is
\begingroup
\allowdisplaybreaks
\begin{align}
\label{homonegx1x2}
f_{X_1,X_2}(x_1,x_2)&\sim \frac{(\tilde{x}_1\tilde{x}_2)^{-\frac{1}{2\eta\xi}-1}}{\sigma_1\sigma_2} \Bigg[\frac{g\left(\omega_x\right)}{4\eta^2}+\left(\tilde{x}_1\tilde{x}_2\right)^{\frac{1}{\xi}}\frac{\tilde{x}_1^{1/\xi}-\tilde{x}_2^{1/\xi}}{\left(\tilde{x}_1^{1/\xi}+\tilde{x}_2^{1/\xi}\right)^3} g^{'}\left(\omega_x\right)-\frac{\left(\tilde{x}_1\tilde{x}_2\right)^{\frac{2}{\xi}}g^{''}\left(\omega_x\right)}{\left(\tilde{x}_1^{1/\xi}+\tilde{x}_2^{1/\xi}\right)^4}\Bigg],
\end{align}
\endgroup
as $x_1\rightarrow x_1^F$ and $x_2 \rightarrow x_2^F$ such that $\omega_x\rightarrow \omega^*_x \in (0,1)$. We now apply the transformation $(X_1,X_2)\rightarrow (R,W)$, where $R=X_1+X_2$ and $W=(\sigma_1+\xi X_1)\{(\sigma_1+\xi X_1)+(\sigma_2+\xi X_2)\}^{-1}$. For $r^F=-(\sigma_1/\xi+\sigma_2/\xi)$ the upper-endpoint of $R$, the density of $(R,W)$ as $r \rightarrow r^F$, is $f_{R,W}(r,w)\sim (-\xi)^{-\frac{1}{\eta\xi}-1}\eta^{-1}(\sigma_1+\sigma_2)^{\frac{1}{\eta\xi}} (r^F-r)^{-\frac{1}{\eta\xi}-1}g^*(w),$ where 
\begin{align*}
g^{*}(w)&=\eta(\sigma_1+\sigma_2)^{\frac{1}{\eta\xi}}(-\xi)^{-1}(\sigma_1\sigma_2)^{\frac{1}{\eta\xi}}\left\{w(1-w)\right\}^{-\frac{1}{2\eta\xi}-1} \Bigg[(4\eta^2)^{-1}g\left(t_w\right)\nonumber\\
&\times\Bigg[(4\eta^2)^{-1}g\left(t_w\right)+\left(\frac{w(1-w)}{\sigma_1\sigma_2}\right)^{\frac{1}{\xi}}\frac{\left(w/\sigma_1\right)^{1/\xi}-\left((1-w)/\sigma_2\right)^{1/\xi}}{\left(\left(w/\sigma_1\right)^{1/\xi}+\left((1-w)/\sigma_2\right)^{1/\xi}\right)^3} g^{'}\left(t_w\right)\\
&-\frac{\left(w(1-w)(\sigma_1\sigma_2)^{-1}\right)^{\frac{2}{\xi}}}{\left(\left(w/\sigma_1\right)^{1/\xi}+\left((1-w)/\sigma_2\right)^{1/\xi}\right)^4}g^{''}\left(t_w\right)\Bigg],
\end{align*} 
and $t_w=\left(w/\sigma_1\right)^{1/\xi}\left\{\left(w/\sigma_1\right)^{1/\xi}+\left((1-w)/\sigma_2\right)^{1/\xi}\right\}^{-1}\in(0,1)$. The support of $W$ is independent of $R$ given that $R$ is above $u=\max\{x_1^F,x_2^F\}$. Consider the survival function of $R$ as $s \rightarrow r^F$, so $s > u$, then
\begin{align}
\label{negK}
\Pr\{ R \geq s\} &\sim \int^\infty_s\int^1_0\frac{(-\xi)^{-\frac{1}{\eta\xi}-1}}{\eta(\sigma_1+\sigma_2)^{-\frac{1}{\eta\xi}}} (r^F-r)^{-\frac{1}{\eta\xi}-1}g^*(w)\mathrm{d}w\mathrm{d}r\sim  K\left(1+\xi \frac{s}{\sigma_1+\sigma_2}\right)^{-\frac{1}{\eta\xi}},
\end{align}
with $K=\int^1_0g^{*}(w)\mathrm{d}w<\infty$.\par
Positive Shape Case ($\xi > 0$). We use the inclusion-exclusion formula to write $\Pr\{R \geq r\}$ in terms of events $\{R \geq r \cap X_1 > u_1\}$ and $\{R \geq r \cap X_2 > u_2\}$ for any fixed constants $u_1,u_2 > 0$.  We first derive $\Pr\{ R \geq r \cap X_1 > u_1\}$ for large $u_1$.
Assume that limit \eqref{heffeq1} holds for $a(y)=\alpha y$ and $b(y)=y^\beta$ for some $\alpha \in [0,1]$ and $\beta \in [0,1]$ for some large $u$. We denote the residual distribution by $G_Z(\cdot)$ and assume it is differentiable with density $g_Z$. Then the joint density of $(Y_1,Y_2)|Y_1 > u$ is ${f_{(Y_1,Y_2)|Y_1 > u}(y_1,y_2)= \exp(-y_1)y_1^{-\beta}g_Z\left((y_2-\alpha y_1)y_1^{-\beta}\right)}$, for $y_1 > u_1$ and $y_2\geq 0$. We now transform to heavy tailed marginals $X_1 \sim \mbox{GPD}(\sigma_1, \xi)$ and $X_2 \sim \mbox{GPD}(\sigma_2,\xi)$ for $\xi > 0$ through the transformation $(Y_1,Y_2)\rightarrow (X_1,X_2)$ where $Y_i=\frac{1}{\xi}\log(1+\xi X_i/\sigma_i)$ for $i \in \{1,2\}$. We also note that $Y_1 > u$ is equivalent to $X_1 > \frac{\sigma_1}{\xi}\left\{\exp(\xi u) - 1\right\}:=u_1$ and so we rewrite the condition as $X_1 > u_1$. The joint density of $(X_1,X_2) |X_1 > u_1$ is ${f_{(X_1,X_2)|X_1 > u_1}(x_1,x_2)=\xi^{\beta}(\sigma_1\sigma_2)^{-1}\tilde{x}_2^{-1}\tilde{x}_1^{-1/\xi-1}\left\{\log\left(\tilde{x}_1\right)\right\}^{-\beta} g_Z\left(z_x^*\right)}$, for large $u$ and where $\tilde{x}_i=1+\xi x_i/\sigma_i$ for $i\in\{1,2\}$ and $z_x^{*}=\xi^{\beta-1}\{\log(\tilde{x}_1)\}^{-\beta}[\log(\tilde{x}_2)-\alpha\log(\tilde{x}_1)]$
where $z_x^* \in \mathbb{R}$ if $\beta < 1$ and $z_x^* \geq -\alpha$, otherwise. A transformation $(X_1,X_2)\rightarrow (R,T=-\log(1-X_1/R)/\log(R))$ gives, for $t \in (0,\infty)$ and as $r\rightarrow \infty$, joint density
\begin{align*}
f_{(R,T)|X_1 > u_1}(r,t)&\sim  r^{-1/\xi-1}r^{1-t}\{\log(r)\}^{1-\beta}\xi^{\beta-1}\xi^{-1/\xi}\sigma_1^{1/\xi}\sigma_2^{-1}\left(1+\xi r^{1-t}/\sigma_2\right)^{-1} g_Z\left(z_t^*\right),
\end{align*}
where $z_t^{*}\sim \xi^{\beta-1}\left[\log\left(1+\xi r^{1-t}/\sigma_2\right)-\alpha\log\left(r\right)\right]\left\{\log\left(r\right)\right\}^{-\beta}$.
Hence, as $r\rightarrow \infty$,
\begin{align}
\label{KG}
f_{R|X_1 > u_1}(r)&\sim  \left(\xi/\sigma_1\right)^{-1/\xi}r^{-1/\xi-1}\int_0^\infty r^{1-t}\{\log(r)\}^{1-\beta}\xi^{\beta-1}\sigma_2^{-1}\left(1+\xi r^{1-t}/\sigma_2\right)^{-1} g_Z\left(z_t^*\right)\mathrm{d}w \nonumber\\
&\sim \sigma_1^{1/\xi}(\xi r)^{-1/\xi-1}\left[\bar{G}_Z\left(z_t^*\right)\right]^\infty_0\sim K_G\sigma_1^{1/\xi}(\xi r)^{-1/\xi-1},
\end{align}
as $r\rightarrow \infty$ and where $K_G$ equals $\bar{G}_Z(0)$ if $\alpha = 0$, $\bar{G}_Z(-\alpha)-\bar{G}_Z(1-\alpha)$ if $\beta = 1$, and $1$ otherwise, and so $\Pr\{R \geq s \cap X_i > u_i\} \sim \exp(-u)K_G\xi^{-1/\xi}\sigma_i^{1/\xi}s^{-1/\xi}$, for $i\in\{1,2\}$, as $s \rightarrow \infty$ and for large $u$. Now consider $\Pr\{R \geq s | (X_1 > u_1 \cap X_2 > u_2)\}$, which we derive using characterisation \eqref{landt}. The density of $(X_1,X_2)$ is given by \eqref{homonegx1x2}. We perform the transformation $(X_1,X_2)\rightarrow (R,W=X_1/R)$ and it follows that $f_{(R,W) |(X_1 > u_1 \cap X_2 > u_2)}(r,w)\sim (\sigma_1\sigma_2)^{-1}rg^*(r,w)$ as $r \rightarrow \infty$ for $w \in [u_1/r,1-u_2/r]$, where
\begin{align}
\label{gstarposshapes}
g^*(r,w)\quad &=\quad\left(1+\xi rw/\sigma_1\right)^{-\frac{1}{2\eta\xi}-1}\left(1+\xi r(1-w)/\sigma_2\right)^{-\frac{1}{2\eta\xi}-1}\nonumber\\
&\times \quad\Bigg[(4\eta^2)^{-1}g(t_{r,w})+t_{r,w}(1-t_{r,w})(2t_{r,w}-1)g^{'}(t_{r,w})-t^2_{r,w}(1-t_{r,w})^2g^{''}(t_{r,w})\Bigg],
\end{align}
and $t_{r,w}=\left(1+\xi r w/\sigma_1\right)^{1/\xi}\left[\left(1+\xi r w/\sigma_1\right)^{1/\xi}+\left(1+\xi r(1-w)/\sigma_2\right)^{1/\xi}\right]^{-1}$. It follows that, as $r\rightarrow \infty$, that $f_{R |(X_1 > u_1 \cap X_2 > u_2)}(r)\sim (\sigma_1\sigma_2)^{-1}rI(r)$,
where $I(r)=\int_{u_1/r}^{1-u_2/r}g^{*}(r,w)\mathrm{d}w$ for $g^*$ in \eqref{gstarposshapes}. To evaluate $I(r)$, we consider two cases, each with $I(r) <\infty$. \par
Case 1. We assume Condition 2. Let $I(r)=I_d(r)+I_1(r)+I_2(r)$, where $I_d(r)=\int^{1-d_2}_{d_1}g^*(r,w)\mathrm{d}w$, $I_1(r)=\int_{u_1/r}^{d_1}g^*(r,w)\mathrm{d}w$, and $I_2=\int^{1-u_2/r}_{1-d_2}g^*(r,w)\mathrm{d}w$, where $d_1$ and $d_2$ are constants chosen such that $d_1 > u_1/r, d_2 > u_2/r$ and $ d_1 < 1-d_2$. We show that, as $r \rightarrow \infty$, we have that $I(r) \sim I_1(r) + I_2(r)$. First, consider $I_d(r)$. As $r \rightarrow \infty$, we have that $t_{r,w}\sim{\left(w/\sigma_1\right)^{1/\xi}\left[{\left(w/\sigma_1\right)^{1/\xi}+\left((1-w)/\sigma_2\right)^{1/\xi}}\right]^{-1}}:=t_w$, and it follows that $g^*(r,w)\sim \xi^{-\frac{1}{\eta\xi}-2}(\sigma_1\sigma_2)^{\frac{1}{2\eta\xi}+1}r^{-\frac{1}{\eta\xi}-2}h_w(w)$, where
\[
h_w(w)=w^{-\frac{1}{2\eta\xi}-1}(1-w)^{-\frac{1}{2\eta\xi}-1}  \left[(4\eta^2)^{-1}g(t_w)+t_{w}(1-t_{w})(2t_{w}-1)g^{'}(t_{w})-t^2_{w}(1-t_{w})^2g^{''}(t_{r,w})\right].
\]
Thus, we have $I_d(r) \sim K_{d} r^{-\frac{1}{\eta\xi}-2}$ for constant $K_d = \xi^{-\frac{1}{\eta\xi}-2}(\sigma_1\sigma_2)^{\frac{1}{2\eta\xi}+1}\int^{1-d_2}_{d_1}h_w(w)\mathrm{d}w>0$. Now, consider $I_1(r)$.  We begin by noting that as $r\rightarrow \infty$ and for $w \in[u_1/r,c_1]$, we have $t_{r,w} \rightarrow 0$. From \eqref{gstarposshapes}, it follows that 
\[
g^*(r,w)\sim K_g\left((4\eta^2)^{-1}-\kappa^2\right)\xi^{-\frac{1}{2\eta\xi}-\frac{\kappa}{\xi}-1}(\sigma_2)^{\frac{1}{2\eta\xi}+\frac{\kappa}{\xi}+1}r^{-\frac{1}{2\eta\xi}-\frac{\kappa}{\xi}-1}\left(1+\xi rw/\sigma_1\right)^{-\frac{1}{2\eta\xi}+\frac{\kappa}{\xi}-1},
\]
for $K_g>0$ defined in Condition 2, and $I_1(r)\sim K_4r^{-\frac{1}{2\eta\xi}-\frac{\kappa}{\xi}-2}$ for 
\[
K_4=\left((2\eta)^{-1}-\kappa\right)^{-1}K_g\left({1}/{4\eta^2}-\kappa^2\right)\xi^{-\frac{1}{2\eta\xi}-\frac{\kappa}{\xi}-1}(\sigma_2)^{\frac{1}{2\eta\xi}+\frac{\kappa}{\xi}+2}\left(1+\xi {u_1}/{\sigma_1}\right)^{-\frac{1}{2\eta\xi}+\frac{\kappa}{\xi}} > 0
\] and where the last line follows as $-\frac{1}{2\eta\xi}+\frac{\kappa}{\xi}<0$. By symmetry, $I_2(r) \sim K_5 r^{-\frac{1}{2\eta\xi}-\frac{\kappa}{\xi}-2}$ for constant 
\[
K_5=\left(\kappa-1/(2\eta)\right)^{-1}K_g\left(1/(4\eta^2)-\kappa^2\right)\xi^{-\frac{1}{2\eta\xi}-\frac{\kappa}{\xi}-1}(\sigma_1)^{\frac{1}{2\eta\xi}+\frac{\kappa}{\xi}+2}\left(1+\xi {u_2}/{\sigma_2}\right)^{-\frac{1}{2\eta\xi}+\frac{\kappa}{\xi}} > 0.
\] 
It follows that $I(r) \sim I_1(r) + I_2(r)\sim (K_4+K_5)r^{-\frac{1}{2\eta\xi}-\frac{\kappa}{\xi}-2}$ as $r\rightarrow \infty$. Hence, 
$\Pr\{R \geq s|(X_1 > u_1) \cap (X_2 > u_2)\}\sim K_6 s^{-\frac{1}{2\eta\xi}-\frac{\kappa}{\xi}}$ as $s\rightarrow \infty$, where $K_6=\xi(\sigma_1\sigma_2)^{-1}(K_4+K_5)(1/2\eta+\kappa)^{-1}>0$.
Recall that $u_i = \frac{\sigma_i}{\xi}\left\{\exp(\xi u) - 1\right\}$ for $i\in\{1,2\}$. From \eqref{landt}, we have $\Pr\{(R \geq r )\cap (X_1 > u_1) \cap (X_2 > u_2)\}\sim K_6\exp\left(-u/\eta\right)g\left(1/2\right)r^{-\frac{1}{2\eta\xi}-\frac{\kappa}{\xi}}$, as $r \rightarrow \infty$ and for large $u$.  Combining all terms, we have $\Pr\{ R \geq r\}\sim K_u^+r^{-1/\xi}$ as $r \rightarrow \infty$, where
\begin{equation}
\label{Kplus}
K^{+}_u=\begin{cases}
\exp(-u)K_G\xi^{-1/\xi}(\sigma_1^{1/\xi}+\sigma_2^{1/\xi})-K_6\exp\left(-u/\eta\right)g\left(1/2\right),\;\;&\text{if}\;\;\frac{1}{2\eta}+\kappa=1,\\
\exp(-u)K_G\xi^{-1/\xi}(\sigma_1^{1/\xi}+\sigma_2^{1/\xi}),\;\;&\text{if}\;\;\frac{1}{2\eta}+\kappa>1,
\end{cases}
\end{equation}
for $K_G$ defined in \eqref{KG}. Note that as $\eta \in [1/2,1]$ and $\kappa < 1/(2\eta)$, we have $\frac{1}{2\eta}+\kappa \geq 1$.\par
Case 2. We now assume that $g$ satisfies Condition 3 with $\eta=1$. Then $I(r)\sim \int^1_0g^*(r,w)\mathrm{d}w$, for $g^*$ defined in \eqref{gstarposshapes} as we have that $g^*(r,w) \sim C_1(r)$ as $w\rightarrow 0$ and $g_2^*(w) \sim C_2(r)$ as $w\rightarrow 1$, where $C_1(r),C_2(r)>0$ are constants with respect to $w$. As $g$ satisfies Condition 3, we have for $s_{r,w}=(1-t_{r,w})/t_{r,w}$ that $g^*(r,w)=\left(1+\xi rw/\sigma_1\right)^{-\frac{1}{\xi}-1}\left(1+\xi r(1-w)/\sigma_2\right)^{-1}s_{r,w} \left\{-H_{12}\left(1,s_{r,w}\right)\right\}$,
and it follows, as $r\rightarrow \infty$,  that
\begin{align*}
I(r)&\sim \int^1_0\left(1+\xi rw/\sigma_1\right)^{-\frac{1}{\xi}-1}\left(1+\xi r(1-w)/\sigma_2\right)^{-1}s_{r,w} \left\{-H_{12}\left(1,s_{r,w}\right)\right\}\mathrm{d}w\\
&\sim \xi^{-\frac{1}{\xi}-1}(\sigma_1\sigma_2)^{\frac{1}{2\xi}}r^{-\frac{1}{\xi}-1}\int^1_0g^*_2(w)\mathrm{d}w,
\end{align*} 
where $t_{r,w}\sim \left(w/\sigma_1\right)^{1/\xi}\left[\left(w/\sigma_1\right)^{1/\xi}+\left((1-w)/\sigma_2\right)^{1/\xi}\right]^{-1}:=t_w$ as $r\rightarrow \infty$ and, for $s_w=(1-t_w)/t_w$, we have
$g^*_2(w)=\left(\sigma_1/\sigma_2\right)^{\frac{1}{2\xi}}\xi^{-1}w^{-\frac{1}{2\xi}-1}(1-w)^{-1}s_w\left\{-H_{12}\left(1,s_w\right)\right\}$. It follows from above and \eqref{landt} that
\begin{equation}
\label{poshomolastterm2}
\Pr\{(R \geq r )\cap (X_1 > u_1) \cap (X_2 > u_2)\}\sim  K\exp(-u)g\left(1/2\right)\xi^{-\frac{1}{\xi}}(\sigma_1\sigma_2)^{\frac{1}{2\xi}}r^{-\frac{1}{\xi}},
\end{equation}
as $r \rightarrow \infty$ and for large $u$ and where $K=\int^1_0g_2^*(w)\mathrm{d}w$.
Combining all terms, we have
\begin{align*}
\Pr\{ R \geq r\}  &\sim\exp(-u)K_G\xi^{-1/\xi}\left(\sigma_1^{1/\xi}+\sigma_2^{1/\xi}\right)r^{-1/\xi}-K\exp(-u)g\left(1/2\right)\xi^{-\frac{1}{\xi}}(\sigma_1\sigma_2)^{\frac{1}{2\xi}}r^{-\frac{1}{\xi}},
\end{align*}
as $r \rightarrow \infty$. Combining Cases 1 and 2 we have that $\Pr\{R \geq s\} \sim K_u^*s^{-\frac{1}{\xi}}$ as $s\rightarrow \infty$, where
\begin{equation}
\label{posKall}
K_u^{*}=K_u^+\text{ for Case 1, }K_u^*=\exp(-u)\xi^{-1/\xi}\left\{K_G(\sigma_1^{1/\xi}+\sigma_2^{1/\xi})-\xi Kg\left(1/2\right)(\sigma_1\sigma_2)^{\frac{1}{2\xi}}\right\}\text{ for Case 2},
\end{equation}
and for $K_u^+$ and $K_G$ defined in \eqref{Kplus} and \eqref{KG}, respectively.\par
\textbf{Proof of Theorem \ref{propzero}.} We derive the joint density of $(X_1,X_2)$ implied by the dependence model given in \eqref{landt} and transform $(X_1,X_2)\rightarrow (R,W)$, where $R=X_1+X_2$ and $W$ is an auxiliary variable. By making assumptions about $g(w)$, we can integrate $W$ out and derive the survival function of $R$. From \eqref{landt} and \eqref{GPDCDF}, we have
  \[
 \Pr\left\{1-\exp\left\{-\frac{X_1}{\sigma_1}\right\} > 1-\frac{1}{x_1}, 1-\exp\left\{-\frac{X_2}{\sigma_2}\right\} > 1-\frac{1}{x_2}\right\} = \frac{\mathcal{L}(x_1+x_2)}{(x_1x_2)^{\frac{1}{2\eta}}}g\left(\frac{x_1}{x_1+x_2}\right),
  \]
 as $x_1,x_2\rightarrow \infty$ such that the limit of $x_1/(x_1+x_2)$ is bounded by $(0,1)$. Under the assumption that $\mathcal{L}(y)$ acts as a constant which can be absorbed by $g$ for $y > v$ for some $v>0$, we have $\Pr\left\{X_1>x_1, X_2>x_2\right\} \sim (\tilde{x}_1\tilde{x}_2)^{-\frac{1}{2\eta}}g\left(\omega_x\right)$ as $x_1,x_2 \rightarrow \infty$ such that $\omega_x = \tilde{x}_1/(\tilde{x}_1+\tilde{x}_2)\rightarrow \omega^*_x \in (0,1)$ and where $\tilde{x}_i=\exp(x_i/\sigma_i)$ for $i\in\{1,2\}$; this implies that $x_2 \sim \sigma_2\left(\frac{x_1}{\sigma_1}+\log\left(\frac{1-\omega^*_x}{\omega^*_x}\right)\right)$ as $x_1 \rightarrow \infty$. Under the assumption that the first and second derivatives of $g$ exist, and the transformation $(X_1,X_2) \rightarrow (R,W)$, the density of $(R,W)$ is ${f_{R,W}(r,w)\sim\{\eta(\sigma_1+\sigma_2)\}^{-1}\exp\left(-r\{\eta(\sigma_1+\sigma_2)\}^{-1}\right)g^{*}(w)},$ for $w\in [-r/\sigma_2,r/\sigma_1]$ and as $r \rightarrow \infty$, and where
\begin{align}
\label{zerogfunc}
g^{*}(w)&=\eta\exp\left(-\frac{(\sigma_2-\sigma_1)w}{2\eta(\sigma_1+\sigma_2)}\right) \Bigg[\frac{g(t_w)}{4\eta^2} +t_w(1-t_w)(2t_w-1) g^{'}\left(t_w\right)-t^2_w(1-t_w)^2g^{''}\left(t_w\right)\Bigg]
\end{align}
and we have $t_w\in(0,1)$, where
\begin{align*}
t_w &=\exp\left(\sigma_2w/(\sigma_1+\sigma_2)\right)\left[\exp\left(\sigma_2w/(\sigma_1+\sigma_2)\right)+\exp\left(-\sigma_1w/(\sigma_1+\sigma_2)\right)\right]^{-1}\\
&=\exp\left(w\right)[\exp\left(w\right)+1]^{-1} \in (0,1),
 \end{align*}
 which follows by multiplying the denominator and numerator of $t_w$ by $\exp(\sigma_1w/(\sigma_1+\sigma_2))$. It follows that with $I(r)=\int^{r/\sigma_1}_{-r/\sigma_2}g^{*}(w)\mathrm{d}w$, as $r\rightarrow \infty$,
\begin{align}
\label{zeroheteromarginalR}
f_{R}(r)=\int^{r/\sigma_1}_{-r/\sigma_2} f_{R,W}(r,w)\mathrm{d}w\sim I(r)(\eta(\sigma_1+\sigma_2))^{-1}\exp\left\{-r(\eta(\sigma_1+\sigma_2))^{-1}\right\}.
\end{align}
 To evaluate $I(r)<\infty$, we make different assumptions on how $g(w)$ behaves; we consider two cases.\par
Case 1. We assume Condition 1. Hence, $I(r)$ in \eqref{zeroheteromarginalR} becomes
\begin{align*}
I(r)&=\int^{r/\sigma_1}_{-r/\sigma_2}\frac{1}{4\eta}\exp\left\{-\frac{(\sigma_2-\sigma_1)w}{2\eta(\sigma_1+\sigma_2)}\right\} \mathrm{d}w\\
&=\begin{cases}
\frac{(\sigma_1+\sigma_2)r}{4\eta\sigma_1\sigma_2},\;\;&\text{if}\;\; \sigma_1 = \sigma_2,\\
\frac{\sigma_1+\sigma_2}{2(\sigma_2-\sigma_1)}\left[\exp\left\{-\frac{(\sigma_1-\sigma_2)r}{2\eta\sigma_2(\sigma_1+\sigma_2)}\right\}-\exp\left\{-\frac{(\sigma_2-\sigma_1)r}{2\eta\sigma_1(\sigma_1+\sigma_2)}\right\}\right],\;\;&\text{if}\;\; \sigma_1 \neq \sigma_2.
\end{cases}
\end{align*}
When $\sigma_1 = \sigma_2 = \sigma$ (say), the marginal density of $R$ is $f_{R}(r)\sim r(4\eta^2\sigma^2)^{-1}\exp\left\{-r(2\eta\sigma)^{-1}\right\}$ as $r \rightarrow \infty$, hence, $\Pr\{ R \geq r\} \sim  r(2\eta\sigma)^{-1}\exp(-r/(2\eta\sigma))$ as $r \rightarrow \infty$. Whereas when $\sigma_1 \neq \sigma_2$, we assume, without loss of generality, that $\sigma_2 > \sigma_1$. Then the marginal density of $R$ is
\begin{align*}
f_{R}(r)&\sim \frac{1}{2\eta(\sigma_2-\sigma_1)}\exp\left\{\frac{(\sigma_2-\sigma_1)r}{2\eta\sigma_2(\sigma_1+\sigma_2)}\right\}\exp\left\{-\frac{r}{\eta(\sigma_1+\sigma_2)}\right\}\sim \frac{1}{2\eta(\sigma_2-\sigma_1)}\exp\left\{-\frac{r}{2\eta\sigma_2}\right\},
\end{align*}
as $r\rightarrow \infty$ and so $\Pr\{ R \geq s\} \sim \sigma_2/(\sigma_2-\sigma_1)\exp\left\{-s/(2\eta\sigma_2)\right\}$ as $s \rightarrow \infty$. By symmetry, this can be written as
$\Pr\{ R \geq r\} \sim [\sigma_{max}/(\sigma_{max}-\sigma_{min})]\exp\left(-r/2\eta\sigma_{max}\right)$,
as $r\rightarrow \infty$ where $\sigma_{max}=\max\{\sigma_1,\sigma_2\}$ and $\sigma_{min}=\min\{\sigma_1,\sigma_2\}$.\par
Case 2. Assume that Conditions 3 and 3a hold with $\eta=1$. Then
\begin{equation}
\label{Isimeq}
I(r)=\int^{r/\sigma_1}_{-r/\sigma_2}g^*(w)\mathrm{d}w\sim \int^\infty_{-\infty}g^*(w)\mathrm{d}w := K,
\end{equation}
for finite constant $K>0$, as $g^*(w)=-\exp(aw) H_{12}\left(1,\exp(-w)\right)$ where $a=\sigma_1/(\sigma_1+\sigma_2)-2\in(-2,-1)$, so $g^*(w)\sim K_{H_1}\exp((a-c_1)w) \rightarrow 0$ as $w\rightarrow \infty$ and $g^*(w)\sim K_{H_2}\exp((a-c_2)w) \rightarrow 0$ as $w\rightarrow -\infty$, where the first limit follows as $a-c_1<-2<0$ and the second follows as $a-c_2>0$.
Hence, \eqref{Isimeq} holds, and it follows that the survival function of $R$ as $s\rightarrow \infty$ is
\begin{align*}
\Pr\{ R \geq s\} \sim \int^\infty_s K(\sigma_1+\sigma_2)^{-1}\exp\left\{-r(\sigma_1+\sigma_2)^{-1}\right\}\mathrm{d}r= K\exp\left\{-s(\sigma_1+\sigma_2)^{-1}\right\}.
\end{align*}\par
 \textbf{Proof of Theorem \ref{prophetero}.} The framework of the proof follows that of the $\xi < 0$ case for Theorem \ref{prophomo}. 
Combining \eqref{landt} and \eqref{GPDCDF}, we have that
  \[
 \Pr\left\{1-\left(1+\xi_1\frac{X_1}{\sigma_1}\right)^{-1/\xi_1} > 1-\frac{1}{x_1}, 1-\left(1+\xi_2\frac{X_2}{\sigma_2}\right)^{-1/\xi_2} > 1-\frac{1}{x_2}\right\} = \frac{\mathcal{L}(x_1+x_2)}{(x_1x_2)^{\frac{1}{2\eta}}}g\left(\frac{x_1}{x_1+x_2}\right),
  \]
  as $x_1,x_2\rightarrow \infty$ such that the limit of $x_1/(x_1+x_2)$ is bounded by $(0,1)$. Under the assumption that $\mathcal{L}(y)$ acts as a constant which can be absorbed by $g$ for $y > v$ for some $v>0$, we have ${\Pr\left\{X_1>x_1, X_2>x_2\right\} \sim \tilde{x}_1^{-\frac{1}{2\eta\xi_1}}\tilde{x}_2^{-\frac{1}{2\eta\xi_2}}g\left(\omega_x\right)}$ for $x_1\rightarrow x_1^F$ and $x_2 \rightarrow x_2^F$, such that $\omega_x= \tilde{x}^{1/\xi_1}_1/(\tilde{x}^{1/\xi_1}_1+\tilde{x}^{1/\xi_2}_2)\rightarrow \omega^*_x \in (0,1)$ and where $\tilde{x}_i=(1+\xi_i x_i/\sigma_i)$ for $i\in\{1,2\}$. Assuming that the first and second derivatives of $g$ exist, and applying the transformation $(X_1,X_2)\rightarrow (R,W)$, where $R=X_1+X_2$ and $W=\left(\sigma_1/\xi_1+X_1\right)\left[\left(\sigma_1/\xi_1+X_1\right) + \left(\sigma_2/\xi_2+X_2\right)\right]^{-1}$, the density of $(R,W)$ as $r\rightarrow r^F$ for $r^F=-(\sigma_1/\xi_1+\sigma_2/\xi_2)$ the upper-endpoint of $R$ is
\begingroup
\allowdisplaybreaks
\begin{align}
\label{frw_xineghetero}
f_{R,W}(r,w) &\sim (-\xi_1)^{-\frac{1}{2\eta\xi_1}-1}(-\xi_2)^{-\frac{1}{2\eta\xi_2}-1}\sigma_1^{\frac{1}{2\eta\xi_1}}\sigma_2^{\frac{1}{2\eta\xi_2}}\left(r^F-r\right)^{-\frac{1}{2\eta\xi_1}-\frac{1}{2\eta\xi_2}-1}w^{-\frac{1}{2\eta\xi_1}-1}(1-w)^{-\frac{1}{2\eta\xi_2}-1}  \nonumber\\
&\times\Bigg[(4\eta^2)^{-1}g\left(t_{r,w}\right)+t_{r,w}(1-t_{r,w})(2t_{r,w}-1)g^{'}\left(t_{r,w}\right)-t_{r,w}^2(1-t_{r,w})^2 g^{''}\left(t_{r,w}\right)\Bigg],
\end{align} 
\endgroup
where $t_{r,w}=\left\{-\xi_1\left(r^F-r\right)w/\sigma_1\right\}^{1/\xi_1}
\left[\left\{-\xi_1\left(r^F-r\right)w/\sigma_1\right\}^{1/\xi_1}+\left\{-\xi_2\left(r^F-r\right)(1-w)/\sigma_2\right\}^{1/\xi_2}\right]^{-1}$.
It follows that the support of $W\in[0,1]$ is independent of $R|R>u$ when $u=\max\{-\sigma_1/\xi_1,-\sigma_2/\xi_2\}$. If $\xi_1 > \xi_2$, we have $t_{r,w} \rightarrow 1$ as $r \rightarrow r^F$, and so require assumptions on how $g(t)$ behaves as $t \rightarrow 1$. We consider two cases:\par
Case 1. Assume that Conditions 2 holds. Then the joint density of $(R,W)$ is
\begin{align*}
f_{R,W}(r,w)\sim K_1K_2\left(r^F-r\right)^{-\frac{1+2\eta\kappa}{2\eta\xi_1}-\frac{1-2\eta\kappa}{2\eta\xi_2}-1}w^{-\frac{1}{2\eta\xi_1}-\frac{\kappa}{\xi_1}-1}(1-w)^{-\frac{1}{2\eta\xi_2}+\frac{\kappa}{\xi_2}-1},
\end{align*}
as $r \rightarrow r^F$, for constants ${K_1=-\left((1+2\eta\kappa)(2\eta\xi_1)^{-1}+(1-2\eta\kappa)(2\eta\xi_2)^{-1}\right)(r^F)^{\frac{1+2\eta\kappa}{2\eta\xi_1}-\frac{1-2\eta\kappa}{2\eta\xi_2}} > 0}$ and $K_2=-K_gK_1^{-1}[(4\eta^2)^{-1}-\kappa^2](-\xi_1)^{-\frac{1}{2\eta\xi_1}-\frac{\kappa}{\xi_1}-1}(-\xi_2)^{-\frac{1}{2\eta\xi_2}+\frac{\kappa}{\xi_2}-1}\sigma_1^{\frac{1}{2\eta\xi_1}\frac{\kappa}{\xi_1}}\sigma_2^{\frac{1}{2\eta\xi_2}-\frac{\kappa}{\xi_2}}>0$. Then
\begin{align}
\label{heteronegK}
\Pr\{ R \geq s\} &\sim  K_{(2)}K_1\int^\infty_s\left(r^F-r\right)^{-\frac{1+2\eta\kappa}{2\eta\xi_1}-\frac{1-2\eta\kappa}{2\eta\xi_2}-1}\mathrm{d}r\sim   K_{(2)}\left(1+\xi_1\xi_2 \frac{s}{\sigma_1\xi_2+\sigma_2\xi_2}\right)^{-\frac{1+2\eta\kappa}{2\eta\xi_1}-\frac{1-2\eta\kappa}{2\eta\xi_2}},
\end{align}
as $s \rightarrow r^F$, and where $ K_{(2)} =K_2B\left(-[(2\eta)^{-1}-\kappa]/\xi_1,-[(2\eta)^{-1}+\kappa]/\xi_2\right)>0$; here $B(\cdot,\cdot)$ denotes the beta function and both of its arguments are positive. The general result follows by replacing $\xi_1$ and $\xi_2$ with $\max\{\xi_1,\xi_2\}$ and $\min\{\xi_1,\xi_2\}$ respectively and using the behaviour of $g$ as $t \rightarrow 0$ as well as $t\rightarrow 1$.\par
Case 2. Assume that Condition 3b holds; recall that $\eta=1$. From \eqref{frw_xineghetero}, the density of $(R,W)$ for $w \in[0,1]$ as $r\rightarrow r^F$ is
\begin{align}
\label{frw_case2_negshapes}
f_{R,W}(r,w)\sim \left(r^F-r\right)^{-\frac{1}{\xi_1}-1}\Bigg[-s_{r,w}H_{12}\left(1,s_{r,w}\right)\Bigg]g^*(w),
\end{align}
where $s_{r,w}=(1-t_{r,w})/t_{r,w}$ and $g^*(w)=(-\xi_1)^{-\frac{1}{\xi_1}-1}(-\xi_2)^{-1}\sigma_1^{\frac{1}{\xi_1}}w^{-\frac{1}{\xi_1}-1}(1-w)^{-1}$. To marginalise $W$ out of \eqref{frw_case2_negshapes}, we transform $(R,W)\rightarrow (V,Z)$, where $V=(r^F-R)/r^F$ and $W=1-V^Z$ for $Z \in (0,\infty)$, and so large $R$ now corresponds to small $V>0$. The joint density of $(V,Z)$ for $z \in (0,\infty)$ as $v\downarrow 0$ is
\begin{align}
f_{V,Z}(v,z) \sim K_2\left(-\log(v)\right)(vr^F)^{-\frac{1}{\xi_1}-\frac{1}{\xi_1}+\frac{1}{\xi_2}-1}v^{z/\xi_2}\bigg[-x_r H_{12}\left(1,x_r(vr^F)^{-1/\xi_1+1/\xi_2}v^{z/\xi_2}\right)\bigg],
\end{align}
where $K_2=(-\xi_1)^{-\frac{1}{\xi_1}-1}(-\xi_2)^{-1}\sigma_1^{\frac{1}{\xi_1}}(r^F)^{-1} > 0$ and $x_i^F=-\xi_i/\sigma_i$ for $i\in\{1,2\}$ and the ratio $x_r := (x_2^F)^{1/\xi_2}/(x_1^F)^{1/\xi_1}$. Consider now the integral
\begin{align*}
I_1(v)&=\frac{x_r}{\xi_2}\int^\infty_0 \log(v)(vr^F)^{-1/\xi_1+1/\xi_2}v^{z/\xi_2}
S_1(v,z)\mathrm{d}z=H_{1}\left(1,\infty\right)-H_{1}\left(1,x_r(vr^F)^{-1/\xi_1+1/\xi_2}\right),
\end{align*}
where $S_1(v,z)=H_{12}\left(1,x_r(vr^F)^{-1/\xi_1+1/\xi_2}v^{z/\xi_2}\right)$. Then as $v \downarrow 0$, 
\begingroup
\allowdisplaybreaks
\begin{align}
f_{V}(v)&\sim K_2\xi_2(vr^F)^{-\frac{1}{\xi_1}-1}I_1(v)\sim K_2\xi_2
H_{1}\left(1,\infty\right)(vr^F)^{-\frac{1}{\xi_1}-1},
\end{align}
\endgroup
which follows as $1/\xi_2 - 1/\xi_1 >0$ and $H_1(1,z) \rightarrow 0$ as $z \rightarrow 0$. Transforming back to $R$, then as $r \rightarrow r^F$,
\[
\Pr\{R \geq r\} \sim K_3\left(1+\xi_1\xi_2 \frac{r}{\sigma_1\xi_2+\sigma_2\xi_2}\right)^{-\frac{1}{\xi_1}}
\]
for constant $K_3= |H_{1}\left(1,\infty\right)|(-\xi_1)^{-\frac{1}{\xi_1}}(r^F)^{-\frac{1}{\xi_1}}\sigma_1^{\frac{1}{\xi_1}}>0$.
The general result follows by a symmetric argument replacing $\xi_1$ and $\xi_2$ with $\xi_{max}=\max\{\xi_1,\xi_2\}$ and $\xi_{min}=\min\{\xi_1,\xi_2\}$ respectively; that is, the distribution function is of the form $\Pr\{R \geq r\} \sim K_{(3b)}\left(1+\xi_1\xi_2 {r}/{(\sigma_1\xi_2+\sigma_2\xi_2)}\right)^{-\frac{1}{\xi_{max}}}$ as $r\rightarrow r^F$ where
\begin{equation}
\label{heteronegK2}
 K_{(3b)}=\begin{cases}
|H_{1}\left(1,\infty\right)|(-\xi_1)^{-\frac{1}{\xi_1}}(r^F)^{-\frac{1}{\xi_1}}\sigma_1^{\frac{1}{\xi_1}}>0,\;\;&\text{if}\;\; \xi_{max} = \xi_1,\\
|H_{2}\left(\infty,1\right)|(-\xi_2)^{-\frac{1}{\xi_2}}(r^F)^{-\frac{1}{\xi_2}}\sigma_2^{\frac{1}{\xi_2}}>0,\;\;&\text{if}\;\; \xi_{max} = \xi_2.
\end{cases}
\end{equation}
\par
\textbf{Proof of Theorem \ref{propdiffsigns}.} 
For $X_1 \sim \mbox{GPD}(\sigma_1,\xi_1)$ and $X_2 \sim \mbox{GPD}(\sigma_2,\xi_2)$ consider four cases:  $(\xi_1 > \xi_2, \xi_2 > 0)$, $(\xi_1 > 0, \xi_2 < 0)$, $(\xi_1 > 0, \xi_2 = 0)$, $(\xi_1 = 0, \xi_2 < 0)$ and $(\xi_1 > \xi_2, \xi_2 > 0)$; the other cases follow by symmetry. We present details only for $\xi_1=0$ and $\xi_2 < 0$ as the other cases are covered by Breiman's lemma \citep{breiman1965some} since in each of the other cases at least one of the shape parameters is positive; so let $X_1 \sim \mbox{GPD}(\sigma_1,0)$ and $X_2 \sim \mbox{GPD}(\sigma_2, \xi_2 <0)$.\par The proof follows by considering the limiting cases of association, i.e., perfect positive dependence and perfect negative dependence, between $X_1$ and $X_2$. The result is then easily extended for any $X_1$ and $X_2$. To illustrate this, let $R_D$ and $R_N$ be $R$ such that $X_1$ and $X_2$ are perfectly positively-dependent, and perfectly negatively-dependent, respectively. For any $y > 0$, we have 
\begin{align}
\label{the4proofeq}
&\min\left\{\Pr\{ R_D \geq y\},\Pr\{ R_N \geq y\}\right\} \leq \Pr\{ R \geq y\} \leq \max\left\{\Pr\{ R_D \geq y\},\Pr\{ R_N \geq y\}\right\}.
\end{align}
We show, as $y\rightarrow \infty$, that $\Pr\{ R_D \geq y\} \sim C_DS(y)$ and $\Pr\{ R_N \geq y\} \sim C_NS(y)$ for some function $S(y)$ and constants $C_D,C_N>0$. Hence by \eqref{the4proofeq}, it follows that $\lim_{y\rightarrow\infty}\Pr\{ R \geq y \}/S(y)=C$ for $C\in[\min\{C_D,C_N\},\max\{C_D,C_N\}]$. 
\par
For perfect positive dependence $X_2 = F^{-1}_2\{F_1(X_1)\}$ and $\Pr\{R_D \geq r\}=\Pr\{X_1+F^{-1}_2\{F_1(X_1)\}\geq r \}=\Pr\{X_1 \geq x^*\}$ where $x^*$ solves $r=x^*+ F^{-1}_2\{F_1(x^*)\}$. We show $\Pr\{X_1 \geq x^*\} \sim \Pr\{X_1 \geq r\}$ as $r\rightarrow \infty$. To solve for $x^*$, we begin with the initial solution $x^*_0=r$ and consider $x^*_1=r+\epsilon$. Using $X_2=\sigma_2\left[-1+\exp\left(X_1/(\sigma_1/\xi_2)\right)\right]/\xi_2$, this gives $\epsilon=-\sigma_2\left[-1+\exp\left(r/(\sigma_1/\xi_2)\right)\right]/\xi_2\rightarrow \sigma_2/\xi_2$ as $r\rightarrow \infty$ as $\xi_2<0$ and hence $x^*\sim r\left\{1+O\left(\exp(r/(\sigma_1\xi_2))/{r}\right)\right\}+\sigma_2/\xi_2$,
as $r\rightarrow \infty$. It follows that, as $r\rightarrow \infty$ and for constant $C_D=\exp(-\sigma_2/(\sigma_1\xi_2))>1$,
\[
\Pr\{R_D \geq r\}= \exp\left(-r\left\{1+O\left(\exp(r/(\sigma_1\xi_2))/{r}\right)\right\}/\sigma_1-\sigma_2/(\sigma_1\xi_2)\right) \sim C_D\Pr\{X_1 \geq r\}.
\]\par
A similar approach is taken for perfect negative dependence; here $X_2 = 1-F^{-1}_2\{F_1(X_1)\}$ and $\Pr\{R_N \geq r\}=\Pr\{X_1+1-F^{-1}_2\{F_1(X_1)\}\geq r \}=\Pr\{X_1 \geq x^*\}$ where $x^*$ solves $r=x^*+ 1-F^{-1}_2\{F_1(x^*)\}$, and we show $\Pr\{X_1 \geq x^*\} \sim \Pr\{X_1 \geq r\}$ as $r\rightarrow \infty$. To solve for $x^*$, we begin with the initial solution $x^*_0=r$ and consider $x^*_1=r+\epsilon$. With $X_2=\sigma_2\left[-1+\left(1-\exp\left(-X_1/\sigma_1\right)\right)^{-1/\xi_2}\right]/\xi_2$, this gives $\epsilon=-\sigma_2\left[-1+\left(1-\exp\left(-r/\sigma_1\right)\right)^{-1/\xi_2}\right]/\xi_2 \rightarrow 0$ as $r\rightarrow \infty$ as $\xi_2<0$ and hence $x^*\sim r\left\{1+O\left(\left[\left(1-\exp(-r/\sigma_1)\right)^{-1/\xi_2}-1\right]/{r}\right)\right\}$,
as $r\rightarrow \infty$. It follows that, as $r\rightarrow \infty$ and for constant $C_N=1$,
\[
\Pr\{R_N \geq r\}= \exp\left(-r\left\{1+O\left(\left[\left(1-\exp(-r/\sigma_1)\right)^{-1/\xi_2}-1\right]/{r}\right)\right\}/\sigma_1\right) \sim C_N\Pr\{X_1 \geq r\}.
\]
Note that as $\xi_2<0$, then $C_D>C_N=1$. Hence by \eqref{the4proofeq}, we have that for general $X_1$ and $X_2$ with any association, that $\Pr\{R > r\}\sim C \Pr\{X_1 \geq r\}\sim C\exp(-r/\sigma_1)$ as $r\rightarrow \infty$ for $C\in[1,C_D^*]$; by symmetry, to also cover the  $X_1\sim \mbox{GPD}(\sigma_1,\xi_1<0)$ and $X_2\sim\mbox{GPD}(\sigma_2,0)$, we have
\begin{equation}
\label{C1eq}
C_D^*=\exp(-\sigma_{min}/(\sigma_{max}\xi_{min}))> 1,
\end{equation}
where $\xi_{min}=\min\{\xi_1,\xi_2\}<0$, $\sigma_{max}=\{\sigma_i; i\text{ is s.t. } \xi_i = 0\}$ and $\sigma_{min}=\{\sigma_i; i\text{ is s.t. } \xi_i = \xi_{min}\}$.
 \section{Linking \eqref{Rforms} to the usual GPD modelling framework}
  \label{Rforms-link-sec}
Assume that \eqref{Rforms} holds in equality, rather than asymptotically (as in $=$ not $\sim$), for $r \geq u_R$ for fixed $u_R \geq 0$. If $\xi_R > 0$, we have $\Pr\{R \geq r\} = K_{1}r^{-1/\xi_R}$ for $r \geq u_R$, and then for $r > 0$
\begin{align*}
\Pr\{R \geq r+u_R | R > u_R\} &= \frac{K_{1}(r+u_R)^{-1/\xi_R}}{K_{1}u_R^{-1/\xi_R}}=\left(1+\frac{r}{u_R}\right)^{-1/\xi_R}=\left(1+\frac{\xi_R r}{u_R\xi_R}\right)^{-1/\xi_R}.
\end{align*}
It follows that $\left(R -u_R\right)|\left( R > u_R\right) \sim \mbox{GPD}(\sigma_R,\xi_R)$, with $\sigma_R=u_R\xi_R$. A similar approach can be used to show that if $\xi_R = 0$, then $\left(R -u_R\right)|\left( R > u_R\right)$ is  $\mbox{GPD}(\sigma_R,0)$. For $\xi_R < 0$ and $r>0$ with $r+u_R < r^F$, we have
\begin{align*}
\Pr\{R \geq r+u_R | R > u_R\} = \frac{K_{3}\left(1-\frac{r+u_R}{r^F}\right)^{-1/\xi_R}}{ K_{3}\left(1-\frac{u_R}{r^F}\right)^{-1/\xi_R}}&=\left(1-\frac{r}{(r^F-u_R)}\right)^{-1/\xi_R}\\
&=\left(1+\xi_R\frac{r}{(-\xi_R)(r^F-u_R)}\right)^{-1/\xi_R},
\end{align*}
and so $\left(R -u_R\right)|\left( R > u_R\right) \sim \mbox{GPD}(\sigma_R,\xi_R)$, with $\sigma_R=(-\xi_R)(r^F-u_R)$. Note that we have made no assumptions about $r^F$ as this is fully determined by the marginal upper-endpoints.
\end{appendix}
\baselineskip=14pt
\begingroup
\setstretch{0.75}
\bibliographystyle{apalike}

\bibliography{ref}

\endgroup

\end{document}